\documentclass[11pt,twoside]{article}
\usepackage{amssymb}
\usepackage{amssymb,amsmath,amsthm,amsfonts,mathrsfs,hyperref}
\usepackage{times}
\usepackage{enumerate}
\usepackage{cite,titletoc}
\allowdisplaybreaks

\def\cQ{\mathcal Q}

\def\cS{\mathcal S}

\def\cW{\mathcal W}
\def\cX{\mathcal X}
\def\cY{\mathcal Y}

\def\N{\mathop{\mathbb N\kern 0pt}\nolimits}
\def\Q{\mathop{\mathbb Q\kern 0pt}\nolimits}
\def\R{\mathop{\mathbb R\kern 0pt}\nolimits}
\def\SS{\mathop{\mathbb S\kern 0pt}\nolimits}

\def\ds{\displaystyle}

\def\p{\partial}
\def\eps{\epsilon}
\def\dl{\delta}
\def\g{\gamma}
\def\ve{\varepsilon}
\def\f{\frac}

\def\la{\lambda}
\def\dive{\operatorname{div}}
\def\curl{\operatorname{curl}}

\def\ls{\lesssim}

\newcommand{\w}[1]{\langle {#1} \rangle}

\hypersetup{colorlinks=true,linkcolor=blue,citecolor=red,urlcolor=cyan}

\theoremstyle{plain}
\newtheorem{theorem}{Theorem}[section]

\newtheorem{lemma}[theorem]{Lemma}
\newtheorem{corollary}[theorem]{Corollary}
\newtheorem*{conjecture}{Conjecture}

\theoremstyle{definition}

\newtheorem{remark}{Remark}[section]

\numberwithin{equation}{section}

\title{Delayed singularity formation for the three dimensional compressible Euler equations with non-zero vorticity}

\author{Fei Hou$^{1, *}$ \qquad Huicheng
  Yin$^{1,2, }$\footnote{Fei Hou (\texttt{fhou$@$nju.edu.cn}) and
    Huicheng Yin (\texttt{huicheng$@$nju.edu.cn}, \texttt{05407$@$njnu.edu.cn}) are supported by
    the NSFC (No.~11731007).}\\
    [12pt] {\small 1. Department of Mathematics, Nanjing University, Nanjing 210093, China}\\
  {\small 2. School of Mathematical Sciences and Mathematical Institute, }\\
  {\small Nanjing Normal University, Nanjing 210023, China}}

\begin{document}

\date{}
\maketitle
\thispagestyle{empty}

\begin{abstract}
For the 3D compressible isentropic Euler equations with an initial perturbation of
size $\ve$ of a rest state, if the initial vorticity is of size
$\dl$ with  $0<\dl\le \ve$ and $\ve$ is small, we establish that the lifespan of the smooth solutions is $T_{\dl}=O(\min\{e^\frac{1}{\ve},\frac{1}{\delta}\})$
for the polytropic gases, and  $T_{\dl}=O(\frac{1}{\delta})$ for the   Chaplygin gases. For example,
when $\dl=e^{-\f{1}{\ve^2}}$ is chosen, then $T_{\dl}=O(e^{\f{1}{\ve}})$ for the polytropic gases
and $T_{\dl}=O(e^{\f{1}{\ve^2}})$ for the   Chaplygin gases although the perturbations of the initial density and
the divergence of the initial velocity are only of order $O(\ve)$.
Our result illustrates that the time of existence of smooth solutions
depends crucially on the size of the vorticity of the initial data,
as long as the initial data is sufficiently close to a constant.
The main ingredients in the paper are: introducing some suitably weighted energies,
deriving the pointwise space-time decay estimates of solutions, looking for the good unknown
instead of the velocity, and establishing the required weighted estimates on the vorticty.

\vskip 0.2 true cm


\vskip 0.2 true cm
\noindent
\end{abstract}

\vskip 0.6 true cm
\tableofcontents

\section{Introduction}
\subsection{Setting of the problem and statement of the main result}
In this paper, we are concerned with the long time existence of smooth solutions to the
3D compressible Euler equations
\begin{equation}\label{Euler}
\left\{
\begin{aligned}
&\p_t\rho+\dive(\rho u)=0,\\
&\p_t(\rho u)+\dive(\rho u \otimes u)+\nabla p=0,\\
&\rho(0,x)=\bar\rho+\rho^0(x),\quad u(0,x)=u^0(x),
\end{aligned}
\right.
\end{equation}
where $(t,x)=(t, x_1, x_2, x_3)\in\R^{1+3}_+:=[0,\infty)\times\R^3$, $\nabla=(\p_{x_1}, \p_{x_2}, \p_{x_3})$, and $u=(u_1,u_2,u_3),~\rho,~p$ stand for the velocity, density, pressure, respectively.
In addition, $\bar\rho>0$ is a constant, $\rho(0,x)>0$, $u^0(x)=(u^0_1(x), u^0_2(x),u^0_3(x))$, and $(\rho^0(x), u^0(x))\in C^\infty(\R^3)$.
Assume that the pressure $p=p(\rho)$ is smooth on its argument $\rho$.

For the polytropic gases (see \cite{CF:book}),
\begin{equation}\label{polytropic:gas}
    p(\rho)=A\rho^\gamma,
\end{equation}
where $A$ and $\gamma$ ($1<\gamma<3$) are some positive constants.

For the Chaplygin gases (see \cite{CF:book} or \cite{Godin07}),
\begin{equation}\label{Chaplygin:gas}
    p(\rho)=P_0-\frac{B}{\rho},
\end{equation}
where $P_0>0$ and $B>0$ are constants.

If $(\rho,u)\in C^1$ and $\rho>0$, then \eqref{Euler} can be reduced to
\begin{equation}\label{EulerC1form}
\left\{
\begin{aligned}
&\p_t\rho+\dive (\rho u)=0,\\
&\p_tu+u\cdot\nabla u+\frac{c^2(\rho)}{\rho}\nabla\rho=0,
\end{aligned}
\right.
\end{equation}
where $c(\rho)=\sqrt{p'(\rho)}$ is the local sound speed. If we let $(\rho, u)(t,x)=(\hat\rho, \hat u)(t,s)$ with $s=x\cdot\omega$, $\omega=(\omega_1, \omega_2, \omega_3)\in\Bbb S^2$, then system \eqref{EulerC1form} becomes
\begin{equation}\label{Euler:planewave}
\left\{
\begin{aligned}
&\p_t\hat\rho+(\omega_1\hat u_1+\omega_2\hat u_2+\omega_3\hat u_3)\p_s\hat\rho
+\hat\rho(\omega_1\p_s\hat u_1+\omega_2\p_s\hat u_2+\omega_3\p_s\hat u_3)=0,\\
&\p_t\hat u_1+(\omega_1\hat u_1+\omega_2\hat u_2+\omega_3\hat u_3)\p_s\hat u_1
+\frac{\omega_1c^2(\hat\rho)}{\hat\rho}\p_s\hat\rho=0,\\
&\p_t\hat u_2+(\omega_1\hat u_1+\omega_2\hat u_2+\omega_3\hat u_3)\p_s\hat u_2
+\frac{\omega_2c^2(\hat\rho)}{\hat\rho}\p_s\hat\rho=0,\\
&\p_t\hat u_3+(\omega_1\hat u_1+\omega_2\hat u_2+\omega_3\hat u_3)\p_s\hat u_3
+\frac{\omega_3c^2(\hat\rho)}{\hat\rho}\p_s\hat\rho=0.
\end{aligned}
\right.
\end{equation}
It follows from direct computation that \eqref{Euler:planewave} has four eigenvalues
\begin{equation*}
    \hat\lambda_1=\omega_1\hat u_1+\omega_2\hat u_2+\omega_3\hat u_3-c(\hat\rho),~~~
    \hat\lambda_{2,3}=\omega_1\hat u_1+\omega_2\hat u_2+\omega_3\hat u_3,~~~
    \hat\lambda_4=\omega_1\hat u_1+\omega_2\hat u_2+\omega_3\hat u_3+c(\hat\rho),
\end{equation*}
and the corresponding four right eigenvectors are
\begin{equation*}
\begin{array}{l}
\ds\hat r_1=(-\hat\rho, c(\hat\rho)\omega_1, c(\hat\rho)\omega_2, c(\hat\rho)\omega_3)^T,\\
\ds\hat r_2=(0, -\omega_2, \omega_1, 0)^T,\\
\ds\hat r_3=(0, -\omega_3, 0, \omega_1)^T,\\
\ds\hat r_4=(\hat\rho, c(\hat\rho)\omega_1, c(\hat\rho)\omega_2, c(\hat\rho)\omega_3)^T.
\end{array}
\end{equation*}
It is easy to verify that for each $\omega\in\Bbb S^2$,
\begin{equation*}
    \nabla_{\hat\rho,\hat u}\hat\lambda_i\cdot \hat r_i\equiv 0,~~i=2,3.
\end{equation*}
For the Chaplygin gases, then
\begin{equation*}
    \nabla_{\hat\rho,\hat u}\hat\lambda_i\cdot \hat r_i\equiv 0,~~i=1,4.
\end{equation*}
For the polytropic gases, then
\begin{equation*}
    \nabla_{\hat\rho,\hat u}\hat\lambda_i\cdot \hat r_i=\hat\rho c'(\hat\rho)+c(\hat\rho)>0
    \quad \text {for $\hat\rho>0$},~~i=1,4.
\end{equation*}
By the definition in Page 89 of \cite{Majda},  \eqref{EulerC1form} is totally linearly degenerate
for the Chaplygin gases, and \eqref{EulerC1form} is genuinely nonlinear with respect to the
first eigenvalue $\la_{1}$ and the fourth eigenvalue $\la_{4}$
for the polytropic gases when $\rho>0$. Generally speaking, for the nonlinear hyperbolic conservation laws with
small initial data or small perturbed initial data, the genuinely nonlinear condition
can arise the blowup of smooth solutions in finite time and corresponds to the formation of shock
(see \cite{Alinhac95}-\cite{Alinhac95-1}, \cite{BSV-1}-\cite{BSV-3}, \cite{Christodoulou07}-
\cite{CM14}, \cite{Godin05}, \cite{Hormander97book}, \cite{Lax}-\cite{LukSpeck18} and \cite{Yin}); for
1D nonlinear hyperbolic conservation laws,
the totally linearly degenerate condition can produce the global smooth small data solutions (see \cite{Li}).
For the multidimensional case of nonlinear hyperbolic conservation laws with totally linearly degenerate condition,
A.~Majda posed the following conjecture on Page 89 of \cite{Majda}:

\begin{conjecture}
If the $d$ ($d\ge 2$) dimensional nonlinear symmetric system is totally linearly degenerate, then it typically has
smooth global solutions when the initial data are in $H^s(\R^d)$ with $s>\frac{d}{2}+1$ unless the solution
itself blows up in finite time.
In particular, the shock wave formation never happens for any smooth initial data.
\end{conjecture}

By our knowledge, so far this conjecture has not been solved yet even for the small initial data.
As illustrated in Page 89 of \cite{Majda}, the above conjecture is mainly of mathematical interest but its resolution would elucidate both the nonlinear nature of the conditions requiring linear degeneracy of each wave field and also might isolate the fashion in which the shock wave formation arises in quasilinear hyperbolic systems. In fact,  with respect to the small perturbed problem of 3D compressible Euler equations
of Chaplygin gases with non-zero vorticity (as the typical model of multidimensional nonlinear symmetric system with totally linearly degenerate condition),
\begin{equation}\label{Euler-Y}
\left\{
\begin{aligned}
&\p_t\rho+\dive(\rho u)=0,\\
&\p_t(\rho u)+\dive(\rho u \otimes u)+\nabla p=0,\\
&\rho(0,x)=\bar\rho+\ve\tilde\rho^0(x),\qquad u(0,x)=\ve\tilde u^0(x),
\end{aligned}
\right.
\end{equation}
it is still completely unknown whether the global smooth solution $(\rho, u)$ of \eqref{Euler-Y} exists or not.

Next we give some illustrations on the irrotational case of \eqref{Euler} with the small perturbed initial data
$(\rho, u)(0,x)=(\bar\rho+\ve\tilde\rho^0(x),\ve\tilde u^0(x)$, where $(\tilde\rho^0(x),\tilde u^0(x))\in C_0^\infty(\R^3)$
are supported in the ball $B(0, R)$. Without loss of generality,
we assume $\bar\rho=c(\bar\rho)=1$. Set the initial vorticity
\begin{equation}\label{curl:def}
\curl \tilde u^0:=(\p_2\tilde u_3^0-\p_3\tilde u_2^0,\p_3\tilde u_1^0-\p_1\tilde u_3^0,\p_1\tilde u_2^0-\p_2\tilde u_1^0)\equiv0,
\end{equation}
then $\curl u(t,x)\equiv0$ always holds as long as $(\rho,u)\in C^1$. This means that
there exists a potential function $\phi$ such that $u=\nabla\phi$.
The Bernoulli's law shows $\p_t\phi+\frac12|\nabla\phi|^2+h(\rho)=0$ with $h'(\rho)=\frac{c^2(\rho)}{\rho}$ and $h(\bar\rho)=0$.
One can easily check that $h(\rho)=\frac{\rho^{\gamma-1}-1}{\gamma-1}$, where  $1<\gamma<3$
for the polytropic gases in \eqref{polytropic:gas}  and $\gamma=-1$ for the Chaplygin gases in \eqref{Chaplygin:gas}.
Then the density $\rho$ can be expressed as
\begin{equation*}
\rho=\big[1-(\gamma-1)(\p_t\phi+\frac12|\nabla\phi|^2)\big]^\frac{1}{\gamma-1}.
\end{equation*}
Substituting this into the first equation of \eqref{Euler} yields
\begin{equation}\label{QLW}
\p_t^2\phi-\triangle\phi+2\sum_{k=1}^3\p_k\phi\p_t\p_k\phi
+(\gamma-1)\p_t\phi\triangle\phi+\sum_{i,j=1}^3\p_i\phi\p_j\phi\p_{ij}^2\phi
+\frac{\gamma-1}{2}|\nabla\phi|^2\triangle\phi=0,
\end{equation}
where the Laplace operator $\ds\Delta:=\sum_{i=1}^3\p_i^2$. In addition, the initial data of $\phi$ are
\begin{equation}\label{QLW:data}
\left\{
\begin{aligned}
&\phi(0,x)=\ve\int_{R}^{x_1}\tilde u_1^0(\eta, x_2, x_3)d\eta,\\
&\p_t\phi(0,x)=-\ve\tilde\rho^0(x)+\ve^2v(x,\ve),
\end{aligned}
\right.
\end{equation}
where $v(x,\ve)=-\frac{1}{2}\ds\sum_{i=1}^3(\tilde u_i^0)^2(x)
-(\tilde\rho^0)^2(x)\int_0^1\Big(\frac{c^2(\rho)}{\rho}\Big)'\bigg|_{\rho=1+\eta\ve\tilde\rho^0(x)}(1-\eta)d\eta$.

For the equation \eqref{QLW} with \eqref{QLW:data}, when $\ve>0$ is small, it follows from Theorem~6.5.3 in \cite{Hormander97book}
that the lifespan $T_\ve$ of smooth solution $\phi$  satisfies $T_\ve\ge e^\frac{C}{\ve}$ with some positive constant $C$.
On the other hand, for the polytropic gases \eqref{polytropic:gas} with $1<\g<3$, it has been known that the lifespan $T_\ve=
O(e^\frac{C}{\ve})$ is optimal (see \cite{Alinhac01b} and \cite{John})
and meanwhile the shock can be formed (see \cite{CM14}, \cite{HKSW}, \cite{Yin}); for the Chaplygin gases  \eqref{Chaplygin:gas} with $\g=-1$,
the lifespan $T_\ve=+\infty$ holds since the corresponding null condition holds (see \cite{Christodoulou86} and \cite{Klainerman}).
In the present paper, we are concerned with such an interesting question:

{\it When $\curl u^0(x)\not\equiv0$ and $\curl u^0(x)=o(\ve)$, what is the lifespan $T_\ve$ of the classical solution
$(\rho,u)$ to the 3D compressible Euler equations with \eqref{polytropic:gas} or \eqref{Chaplygin:gas}?}

So far, by the author's knowledge, only a few results on the above problem have been obtained for the
2D or 3D compressible Euler equations \eqref{Euler} with non-zero vorticity. For examples, with respect to
the 2D compressible Euler equations of polytropic gases with the following rotationally invariant initial data
\begin{equation}\label{intial:sym}
\rho(0,x)=\bar\rho+\ve\rho^0(r),
\quad u(0,x)=\ve u_r^0(r)\frac{x}{r}+\ve u_\theta^0(r)\frac{x^\perp}{r},
\end{equation}
where $r=|x|$ and $x^\perp=(-x^2,x^1)$, \cite{Alinhac93} has shown that the lifespan $T_\ve$ is of order $O(\frac{1}{\ve^2})$.
The global existence of smooth solution to 2D compressible Euler equations \eqref{Euler} of Chaplygin gases with the initial data \eqref{intial:sym} was established in \cite{HouYin19,HouYin20}. Without the assumption on the rotationally invariant initial data, if $\curl u^0(x)=O(\ve^{1+\alpha})$ with the constant $\alpha\ge0$, it follows from Theorem~1 and Theorem~2 of \cite{Sideris97} that the lifespan $T_\ve$ of 2D compressible Euler equations \eqref{Euler} with \eqref{polytropic:gas} fulfills $T_{\ve}=O(\frac{1}{\ve^{\min\{1+\alpha,2\}}})$.
In addition, for the 3D compressible Euler equations \eqref{EulerC1form} of polytropic gases with the small initial perturbed density of
order $O(\ve)$, when the  divergence of initial velocity is of $O(1)$, and the vorticity of initial velocity is of order $O(\ve^\mu)$
with $1<\mu<\f65$, the authors in \cite{Secchi04} proved that the lifespan of smooth solution $(\rho,u)$ is of $O(\f{1}{\ve^\mu})$.

We now investigate the influence of the vorticity on the  existence time of smooth solutions to 3D compressible Euler equations with \eqref{polytropic:gas} or \eqref{Chaplygin:gas}. To this end, we introduce the following two quantities that capture the size of the perturbed initial data $(\rho^0,u^0)$ and initial vorticity $\curl u^0$:
\begin{equation}\label{initial}
\begin{split}
\ve:=&\sum_{k\le N}\|(\w{|x|}\nabla)^k(\rho^0(x), u^0(x))\|_{L^2},\\
\delta:=&\sum_{k\le N-1}\|\w{|x|}(\w{|x|}\nabla)^k\curl u^0(x)\|_{L^2},
\end{split}
\end{equation}
where $N\ge8$, and $\w{|x|}=\sqrt{1+|x|^2}$. The main result in the paper is

\begin{theorem}\label{thm:main}
For the numbers $\ve$ and $\dl$ defined in \eqref{initial}, then there exist three constants $\ve_0,\delta_0,\kappa_0>0$ such that
when $\ve\le\ve_0$ and $\delta\le\delta_0$,

(i) \eqref{Euler} with the state equation \eqref{polytropic:gas} of polytropic gases  has a solution $(\rho-\bar\rho,u)\in C([0,T], H^N(\R^3))$, where $T=\min\{e^\frac{\kappa_0}{\ve}-1,\frac{\kappa_0}{\delta}\}$.

(ii) \eqref{Euler}  with the  state equation  \eqref{Chaplygin:gas} of Chaplygin gases has a solution $(\rho-\bar\rho,u)\in C([0,T], H^N(\R^3))$, where $T=\frac{\kappa_0}{\delta}$.
\end{theorem}

\subsection{Remarks and sketch of proof}

\begin{remark}
For the Euler-Maxwell system, the authors in \cite{GIP16} have shown that the smooth, irrotational,
small-amplitude solution to the Euler-Maxwell two-fluid system globally exists.
When the initial data are of order $O(\ve)$  and the initial vorticity is of order $O(\dl)$ with
$0<\delta\le\ve$, Ionescu and Lie in \cite{IL18} have proved that the existence time is larger than
$\frac{C}{\delta}$, where $C$ is some positive constant.
\end{remark}

\begin{remark}
When $\delta=e^{-\ve^\ell}$ with $\ell\ge1$ or $e^{-e^\frac{1}{\ve^p}}$ or $p>0$ are chosen in Theorem~\ref{thm:main},
we know that the existence time of smooth solution $(\rho,u)$ to \eqref{Euler} for the Chaplygin gases is larger than
$\kappa_0e^\frac{1}{\ve^\ell}$ or $\kappa_0e^{e^\frac{1}{\ve^p}}$, which means that
the order of lifespan $T_{\delta}$ of small perturbed solution $(\rho, u)$ is only essentially influenced by the
size of the initial vorticity.
\end{remark}

\begin{remark}
In \cite{Secchi04}, through decomposing the solution $(\rho, u)$ of 3D compressible Euler equations \eqref{EulerC1form}
as the sum of the irrotational part, the incompressible part and the remainder and by applying the direct  energy method,
the authors establish that the lifespan of smooth solution $(\rho,u)$ is of order $O(\frac{1}{\ve^\mu})$ ($1<\mu<\frac65$)
when the initial vorticity is of $O(\ve^\mu)$ and the perturbed density is of $O(\ve)$.
We point out that our ingredients in the paper are different from the ones in \cite{Secchi04}.
\end{remark}

\begin{remark}
In our paper \cite{HouYin21}, for the 2D compressible isentropic Euler equations \eqref{Euler} of Chaplygin gases, when
the initial data are a perturbation of size $\ve$, and the initial vorticity is of any size
$\delta$ with  $0<\delta\le \ve$, we have established  the lifespan $T_{\delta}=O(\frac{1}{\delta})$. The main methods in \cite{HouYin21}
are: establishing a new class of  weighted space-time $L^\infty$-$L^\infty$ estimates for the solution itself
and its gradients of 2D linear wave equations,  introducing some suitably weighted energies
and taking the $L^p$ $(1<p<\infty)$ estimates on the vorticity due to the requirements of Sobolev embedding theorem.
Since the space-time decay rates of smooth solutions to 3D and 2D wave equations are different, moreover,
the vorticity equations in 3D and 2D Euler equations have also some differences (for examples, the vorticity $\curl u$ in 2D Euler equations
is a scalar function and there exists the conservation form $(\p_t+u\cdot\nabla)(\frac{\curl u}{\rho})\equiv 0$.
However, for the 3D Euler equations, the vorticity $\curl u$ is a 3D vector and satisfies
a nonlinear system $(\p_t+u\cdot\nabla)\curl u
=\curl u\cdot\nabla u-\curl u\dive u$), then there are a few different techniques and methods between the
present paper and
 \cite{HouYin21} in order to derive the lifespan of smooth solutions.
\end{remark}

Now we give some illustrations  on the proof of Theorem~\ref{thm:main}.
At first, we introduce the perturbed sound speed $\sigma=\frac{c(\rho)-1 }{\lambda}$ as the new unknown
to rewrite  \eqref{Euler} as
\begin{equation}\label{reducedEuler}
\left\{
\begin{aligned}
&\p_t\sigma+\dive u=Q_1:=-\lambda\sigma\dive u-u\cdot\nabla\sigma,\\
&\p_tu+\nabla\sigma=Q_2:=-\lambda\sigma\nabla\sigma-u\cdot\nabla u,
\end{aligned}
\right.
\end{equation}
where $\lambda=\frac{\gamma-1}{2}$, and the $i$-component of the vector $Q_2$ is $Q_{2i}=-\lambda\sigma\p_i\sigma-u\cdot\nabla u_i$.
Secondly, as in \cite{HouYin21}, we introduce the good unknown $g$ in the region $|x|>0$,
\begin{equation}\label{goodunknown:def}
g:=(g_1,g_2,g_3)=u-\omega\sigma\quad\text{with $g_i=u_i-\sigma\omega_i,\quad i=1,2,3,$}
\end{equation}
where $\omega=(\frac{x_1}{|x|},\frac{x_2}{|x|},\frac{x_3}{|x|})\in\SS^2$. Note that the introduction of
$g$ is motivated by the second order quasilinear wave equation \eqref{QLW} although  \eqref{Euler} admits the non-zero vorticity
and can not be transformed a wave equation directly: in  \eqref{QLW},
due to $u_i=\p_i\phi$ and $\sigma=-\p_t\phi+\text{higher order error terms of $\p\phi$}$,
then  $g_i=(\p_i+\omega_i\p_t)\phi+\text{higher order error terms of $\p\phi$}$.
It is well known that $(\p_i+\omega_i\p_t)\phi$ is the good derivative in the study of the nonlinear wave
equation (see \cite{Alinhac10}) since $(\p_i+\omega_i\p_t)\phi$ will admit more rapid
space-time decay rates. By some ideas and methods dealing with the null condition structures
in \cite{HouYin19,HouYin20,HouYin20jde} for the non-compactly supported solutions of
2D quasilinear wave equations, we can obtain better $L^\infty$ space-time decay rates of $g$.
On the other hand, since the optimal time-decay rate
of solutions to the 3D free wave equation is merely $(1+t)^{-1}$ near the forward light cone surface,
which is far to derive the existence time $\ds T_\delta
=\frac{\kappa}{\delta}$ in Theorem  \ref{thm:main} (ii). In fact, when such a $\delta=e^{-e^{\frac{1}{\ve^2}}}$ is chosen,
then the integral $\ds\int_0^{T_{\delta}}\frac{dt}{1+t}=O(e^{\frac{1}{\ve^2}})$ is sufficiently large  as $\ve\rightarrow0$,
which leads to that the usual energy $E(t)$ can not be controlled well by the corresponding
inequality $E(t)\le E(0)+\f{C\ve}{1+t}E(t)$. To overcome this difficulty, some better space-time weighted $L^\infty$-$L^\infty$ estimates
of $g$ and $(\sigma, u)$ are obtained by introducing some auxiliary energies
(including the space-time weighted energies of $\dive u$ and $\curl u$) and looking for the null condition in the system \eqref{Euler}
for the Chaplygin gases, meanwhile the suitable energy estimates of the vorticity are also derived. Based on these key estimates, Theorem ~\ref{thm:main} can be eventually proved.

This paper is organized as follows.
In Section~\ref{sect2}, we will introduce the basic bootstrap assumptions, Helmholtz decomposition and some pointwise estimates.
The estimates of the auxiliary energies and the good unknown $g$ are established in Section~\ref{sect3}.
Collecting the pointwise space-time estimates in Section~\ref{sect2} and \ref{sect3}, the Hardy inequality and the ghost weight method in \cite{Alinhac01a}, we derive the related energy estimates in Section~\ref{sect4}.
In Section~\ref{sect5}, based on the previous energy inequalities and Gronwall's inequalities,
the proof of Theorem~\ref{thm:main} is finished by the continuity argument.

\section{Some preliminaries}\label{sect2}

\subsection{The vector fields and bootstrap assumptions}
Define the spatial rotation vector fields
\begin{equation*}
\Omega:=x\wedge\nabla=(\Omega_{23},\Omega_{31},\Omega_{12}),
\quad \Omega_{ij}:=x_i\p_j-x_j\p_i.
\end{equation*}
For a 3D vector-valued function $U$, denote
\begin{equation*}
\tilde\Omega_{ij}U:=\Omega_{ij}U+e_i\otimes e_jU-e_j\otimes e_iU,
\end{equation*}
where $e_i:=(0,\cdots,\stackrel{i}{1},\cdots,0)^T$.
Define $\tilde\Omega=\{\tilde\Omega_{ij}\}$.
Let $\tilde\Omega U_i=(\tilde\Omega U)_i$ be the $i$-component of $\tilde\Omega U$ rather than the operator $\tilde\Omega$ acts
on the component $U_i$.

According to the definitions of $\Omega$ and $\tilde\Omega$, it is easy to check that for the scalar function $f$ and the
3D vector-valued functions $U,V$,
\begin{equation}\label{rotation:commutation}
\begin{array}{ll}
  \Omega\dive U=\dive\tilde\Omega U,\qquad
  & \tilde\Omega\curl U=\curl\tilde\Omega U, \\
  \tilde\Omega\nabla f=\nabla\Omega f,\qquad
  & \tilde\Omega(U\cdot\nabla V)=U\cdot\nabla(\tilde\Omega V)+(\tilde\Omega U)\cdot\nabla V.
\end{array}
\end{equation}
The spatial derivatives can be decomposed into the radial and angular components for $r=|x|\neq0$,
\begin{equation*}
\nabla=\omega\p_r-\frac{1}{|x|}\omega\wedge\Omega.
\end{equation*}
For convenience, we schematically denote this decomposition as
\begin{equation}\label{radial:angular}
\p_i=\omega_i\p_r+\frac{1}{|x|}\Omega.
\end{equation}
For the multi-index $a$, let
\begin{equation}\label{vectorfield}
\cS:=t\p_t+r\p_r,\quad \Gamma^a=\cS^{a_s}Z^{a_z},\quad Z\in\{\p_t,\nabla,\Omega\},
\qquad\tilde\Gamma^a=\cS^{a_s}\tilde Z^{a_z},
\quad \tilde Z\in\{\p_t,\nabla,\tilde\Omega\}.
\end{equation}
As a consequence, by acting $(\cS+1)^{a_s}Z^{a_z}$ on the first equation and $(\cS+1)^{a_s}\tilde Z^{a_z}$ on
the second equation in \eqref{reducedEuler},
we can find the equations of $(\Gamma^a\sigma,\tilde\Gamma^au)$ as follows
\begin{equation}\label{high:eqn}
\begin{split}
\left\{
\begin{aligned}
&\p_t\Gamma^a\sigma+\dive\tilde\Gamma^au=\cQ_1^a:=\sum_{b+c=a}C^a_{bc}Q_1^{bc},\\
&\p_t\tilde\Gamma^au+\nabla\Gamma^a\sigma=\cQ_2^a:=\sum_{b+c=a}C^a_{bc}Q_2^{bc},
\end{aligned}
\right.
\end{split}
\end{equation}
where $C^a_{bc}$ are some suitable constants ($C^a_{a0}=C^a_{0a}=1$) and
\begin{equation}\label{Qbc:def}
\begin{split}
Q_1^{bc}:=&-\lambda\Gamma^b\sigma\dive\tilde\Gamma^cu
-\tilde\Gamma^bu\cdot\nabla\Gamma^c\sigma,\\
Q_2^{bc}:=&-\lambda\Gamma^b\sigma\nabla\Gamma^c\sigma
-\tilde\Gamma^bu\cdot\nabla\tilde\Gamma^cu.
\end{split}
\end{equation}
For integers $m,m_1,m_2\in\N$ with $m_1\ge1$ and $m_2\ge 2$, set
\begin{equation}\label{energy:def}
\begin{split}
E_m(t):=&~\sum_{|a|\le m}\|(\tilde\Gamma^au,\Gamma^a\sigma)(t,x)\|_{L_x^2},\\
\cX_{m_1}(t):=&~\sum_{|a|\le m_1-1}
\|\w{|x|-t}(\dive\tilde\Gamma^au,\p_t\tilde\Gamma^au,\nabla\Gamma^a\sigma,
\p_t\Gamma^a\sigma)(t,x)\|_{L_x^2},\\
\cY_{m_2}(t):=&~\sum_{|a|\le m_2-2}
\|\w{t}^2(\nabla\dive\tilde\Gamma^au,\nabla^2\Gamma^a\sigma)(t,x)\|_{L^2(|x|\le\w{t}/2)},\\
\cW_m(t):=&~\sum_{|a|\le m}\|\w{|x|}\curl\tilde\Gamma^au(t,x)\|_{L_x^2}.
\end{split}
\end{equation}

Choose the integer $N_1$ such that $N_1+3\le N\le2N_1-2$.
Throughout the whole paper, we make the following bootstrap assumptions: for $t\delta\le\kappa_0$,
\begin{equation}\label{bootstrap}
\begin{split}
&E_N(t)+\cX_N(t)\le M\ve,\qquad\cY_N(t)\le M\ve+M\delta(1+t)^{1+M'\ve},\\
&\cW_{N_1}(t)\le M\delta,\qquad\qquad\cW_{N-1}(t)\le M\delta(1+t)^{M'\ve},\\
&\delta\le\ve,\qquad 0<M'\ve\le\frac18,
\qquad M(\ve+\kappa_0)\le1,
\end{split}
\end{equation}
where the constants  $M\ge 1$, $M'>0$ and $\kappa_0>0$ will be chosen. In Section~\ref{sect5}, we will prove that the constant $M$ on the
right hands of the first two lines in \eqref{bootstrap}
can be improved to $\frac12 M$.

\subsection{The Helmholtz decomposition and commutators}
For the 3D rapidly decaying vector function $U=(U_1,U_2,U_3)$, we divide it into the curl-free
part $P_1U$ (irrotational) and the divergence-free part $P_2U$ (solenoidal), which is called Helmholtz decomposition
\begin{equation}\label{Helmholtz}
U=P_1U+P_2U:=-\nabla(-\Delta)^{-1}\dive U+(-\Delta)^{-1}\curl^2U.
\end{equation}
Note that
\begin{equation}\label{double:curl}
\curl^2U=\curl^2P_2U=-\Delta P_2U+\nabla\dive P_2U=-\Delta P_2U.
\end{equation}
We now give some useful inequalities.
\begin{lemma} For the 3D vector function $U$, it holds that
\begin{equation}\label{div:curl:ineq}
\begin{split}
\|\nabla U\|_{L_x^2}&\ls\|\dive U\|_{L_x^2}+\|\curl U\|_{L_x^2},\\
\|\w{|x|-t}\nabla U\|_{L_x^2}&\ls\|\w{|x|-t}\dive U\|_{L_x^2}
+\w{t}\|\w{|x|}\curl U\|_{L_x^2}+\|U\|_{L_x^2}.
\end{split}
\end{equation}
\end{lemma}
\begin{proof}
Since the proof of \eqref{div:curl:ineq} follows from the integration by parts directly, we omit it here.
\end{proof}

\begin{lemma}
For the vector fields $\tilde\Gamma$ defined in \eqref{vectorfield}, we have
the commutators $[\tilde\Gamma,P_1]:=\tilde\Gamma P_1-P_1\tilde\Gamma=0$ and $[\tilde\Gamma,P_2]=0$.
\end{lemma}
\begin{proof} We only prove $\tilde\Gamma P_2U=P_2(\tilde\Gamma U)$ since it is easy to know
$[\tilde\Gamma,P_1]=[\tilde\Gamma,{\rm Id}-P_2]=0$ if $[\tilde\Gamma,P_2]=0$.

For $\tilde\Gamma\in\{\p_t,\nabla\}$, then $[\tilde\Gamma,P_2]=0$ is obvious. We now focus on the case of $\tilde\Gamma\in\{\tilde\Omega,\cS\}$.

According to \eqref{rotation:commutation} and \eqref{double:curl}, we have
\begin{equation*}
-\Delta P_2(\tilde\Omega U)=\curl^2\tilde\Omega U=\tilde\Omega\curl^2 U
=\tilde\Omega(-\Delta)P_2U=-\Delta(\tilde\Omega P_2U).
\end{equation*}
This, together with the uniqueness of the solution to the equation $\Delta w=0$ when $w$ suitably decays,
yields $P_2(\tilde\Omega U)=\tilde\Omega P_2U$.

Analogously, $P_2(\cS U)=\cS P_2U$ comes from
\begin{equation*}
-\Delta P_2(\cS U)=\curl^2\cS U=(\cS+2)\curl^2 U=(\cS+2)(-\Delta)P_2U=-\Delta(\cS P_2U).
\end{equation*}
\end{proof}

\subsection{The pointwise estimates}
\begin{lemma}\label{lem:pointwise}
For the multi-indices $a,b$ with $|a|\le N-2$ and $b\le N-3$, it holds that
\begin{equation}\label{pointwise:curl}
\w{|x|}^2|\curl\tilde\Gamma^au(t,x)|\ls\cW_{|a|+2}(t),
\end{equation}
and
\begin{equation}\label{pointwise:wave}
\begin{split}
\w{|x|}\w{|x|-t}^\frac12(|\tilde\Gamma^au(t,x)|+|\Gamma^a\sigma(t,x)|)\ls
&~E_{|a|+2}(t)+\cX_{|a|+2}(t)+\w{t}\cW_{|a|+1}(t),\\
\w{|x|}\w{|x|-t}(|\nabla\tilde\Gamma^bu(t,x)|+|\nabla\Gamma^b\sigma(t,x)|)\ls
&~E_{|b|+3}(t)+\cX_{|b|+3}(t)+\w{t}\cW_{|b|+2}(t).
\end{split}
\end{equation}
\end{lemma}
\begin{proof}
The proofs of \eqref{pointwise:curl} and \eqref{pointwise:wave} are motivated by Lemma~3.3 in \cite{Sideris00}.
Note that the main difference between \eqref{pointwise:wave} and Proposition~3.3 in \cite{Sideris00}
lies in the appearance of the vorticity on the right hand side.

Recall the Sobolev-type inequalities (3.14b), (3.14c) and (3.14d) in \cite{Sideris00} that for any $|x|>0$,
\begin{equation}\label{Sobolev:ineq1}
\begin{split}
&|x||W(t,x)|\ls
\|\p_r\tilde\Omega^{\le1}W(t,y)\|_{L^2(|y|\ge|x|)}
+\|\tilde\Omega^{\le2}W(t,y)\|_{L^2(|y|\ge|x|)},\\
&|x|\w{|x|-t}^\frac12|U(t,x)|\ls
\|\w{|x|-t}\p_r\tilde\Omega^{\le1}U(t,y)\|_{L^2(|y|\ge|x|)}
+\|\tilde\Omega^{\le2}U(t,y)\|_{L^2(|y|\ge|x|)},\\
&|x|\w{|x|-t}|V(t,x)|\ls
\|\w{|x|-t}\p_r\tilde\Omega^{\le1}V(t,y)\|_{L^2(|y|\ge|x|)}
+\|\w{|x|-t}\tilde\Omega^{\le2}V(t,y)\|_{L^2(|y|\ge|x|)},
\end{split}
\end{equation}
where $\tilde\Omega^{\le m}:=\sum_{0\le|a|\le m}\tilde\Omega^a$ and $W,U,V$ can be 3D vectors
or scalar functions (for scalar functions, $\tilde\Omega$ in \eqref{Sobolev:ineq1} is replaced by $\Omega$).

At first, we deal with \eqref{pointwise:curl} in the region $|x|\ge1/4$.
Let $W=\tilde\Gamma^a\curl u$ in the first inequality of \eqref{Sobolev:ineq1}. Then
\begin{equation*}
\begin{split}
&\;\quad|x|^2|\curl\tilde\Gamma^au(t,x)|\\
&\ls\|\w{|y|}\p_r\tilde\Omega^{\le1}\curl\tilde\Gamma^au(t,y)\|_{L^2(|y|\ge|x|)}
+\|\w{|y|}\tilde\Omega^{\le2}\curl\tilde\Gamma^au(t,y)\|_{L^2(|y|\ge|x|)}\\
&\ls\cW_{|a|+2}(t),
\end{split}
\end{equation*}
which derives \eqref{pointwise:curl} for $|x|\ge1/4$.
On the other hand, in the domain of $|x|\le1/4$, \eqref{pointwise:curl} follows from the standard Sobolev embedding
theorem directly.

The proof of \eqref{pointwise:wave} in the region $|x|\ge1/4$ follows from the choices of $U=\tilde\Gamma^au,\Gamma^a\sigma$, $V=\nabla\tilde\Gamma^au,\nabla\Gamma^a\sigma$ in \eqref{Sobolev:ineq1} and the weighted inequality \eqref{div:curl:ineq}.

Finally, we deal with \eqref{pointwise:wave} in the domain of $|x|\le1/4$.
Choosing the cut-off function $\chi(s)\in C^\infty$ such that
\begin{equation}\label{cutoff}
\begin{split}
0\le\chi(s)\le1,\qquad \chi(s)=\left\{
\begin{aligned}
&1,\qquad\quad s\le1/4,\\
&0,\qquad\quad s\ge1/2.
\end{aligned}
\right.
\end{split}
\end{equation}
Applying the standard Sobolev embedding theorem to $\chi(|x|)\tilde\Gamma^au(t,x)$, one has
\begin{equation}\label{Sobolev:ineq2}
|\chi(|x|)\tilde\Gamma^au(t,x)|
\ls\|\tilde\Gamma^au(t,x)\|_{L^2(|x|\le1)}
+\w{t}^{-1}\|\w{|x|-t}\nabla\nabla^{\le1}\tilde\Gamma^au(t,x)|\|_{L^2(|x|\le1)}.
\end{equation}
By using \eqref{Sobolev:ineq1} to the first term on the right hand side of \eqref{Sobolev:ineq2}, we arrive at
\begin{equation}\label{Sobolev:ineq3}
\begin{split}
\w{t}^\frac12\|\tilde\Gamma^au(t,x)\|_{L^2(|x|\le1)}
&\ls\Big\||x|^{-1}\Big\|_{L^2(|x|\le1)}
\Big\||x|\w{|x|-t}^\frac12\tilde\Gamma^au(t,x)\Big\|_{L^\infty(0<|x|\le1)}\\
&\ls\|\w{|x|-t}\p_r\tilde\Omega^{\le1}\tilde\Gamma^au(t,y)\|_{L^2}
+\|\tilde\Omega^{\le2}\tilde\Gamma^au(t,y)\|_{L^2}.
\end{split}
\end{equation}
Substituting \eqref{Sobolev:ineq3} into \eqref{Sobolev:ineq2} yields  \eqref{pointwise:wave} for $\tilde\Gamma^au$.
On the other hand, for $\Gamma^a\sigma,\nabla\tilde\Gamma^au$ and $\nabla\Gamma^a\sigma$,
\eqref{pointwise:wave} can be analogously proved.
\end{proof}

\begin{lemma}[Sharp time decay of $P_1u$ away from the conic surface $|x|=\w{t}$]\label{lem:sharp:decay}
For $|a|\le N-2$ and $|x|\le\w{t}/4$, it holds that
\begin{equation}\label{sharp:decay1}
\w{t}^\frac32(|P_1\tilde\Gamma^au(t,x)|+|\Gamma^a\sigma(t,x)|)\ls
E_{|a|+2}(t)+\cX_{|a|+2}(t)+\cY_{|a|+2}(t)+\w{t}\cW_{|a|}(t).
\end{equation}
\end{lemma}
\begin{proof}
At first, we consider the general function $\chi(\frac{|x|}{\w{t}})U(t,x)$, where the cut-off
function $\chi$ is defined in \eqref{cutoff}.
Let $x=\w{t}y$ and then it follows from the standard Sobolev embedding theorem
$H_y^2(\R^3)\hookrightarrow L_y^\infty(\R^3)$ that
\begin{equation}\label{sharp:decay2}
\begin{split}
&\;\quad\|\chi(\frac{|x|}{\w{t}})U(t,x)\|_{L_x^\infty}=\|\chi(|y|)U(t,\w{t}y)\|_{L_y^\infty}\\
&\ls\|U(t,\w{t}y)\|_{L_y^2}+\|\w{t}(\nabla_xU)(t,\w{t}y)\|_{L^2(|y|\le1/2)}
+\|\w{t}^2\chi(|y|)(\nabla_x^2U)(t,\w{t}y)\|_{L_y^2}\\
&\ls\w{t}^{-\frac32}\Big\{\|U(t,x)\|_{L_x^2}
+\|\w{t}\nabla U(t,x)\|_{L^2(|x|\le\w{t}/2)}
+\w{t}^2\Big\|\chi(\frac{|x|}{\w{t}})\nabla^2U(t,x)\Big\|_{L_x^2}\Big\}.
\end{split}
\end{equation}
Here we point out that $U(t,x)=P_1\tilde\Gamma^au(t,x)$ or $U(t,x)=\Gamma^a\sigma(t,x)$ will be chosen in \eqref{sharp:decay2}.
Note that the second term in the last line of \eqref{sharp:decay2} can be estimated as follows
\begin{equation}\label{sharp:decay3}
\begin{split}
\|\w{t}\nabla P_1\tilde\Gamma^au\|_{L^2(|x|\le\w{t}/2)}
&\ls\|\w{|x|-t}\nabla\tilde\Gamma^au\|_{L_x^2}
+\w{t}\|\nabla P_2\tilde\Gamma^au\|_{L_x^2}\\
&\ls\cX_{|a|+1}(t)+\w{t}\cW_{|a|}(t)+E_{|a|}(t)+\w{t}\|\curl P_2\tilde\Gamma^au\|_{L_x^2}\\
&\ls E_{|a|+2}(t)+\cX_{|a|+2}(t)+\w{t}\cW_{|a|}(t),
\end{split}
\end{equation}
where we have used the second inequality in \eqref{div:curl:ineq}.
We next deal with the third term in the last line of \eqref{sharp:decay2} with $U(t,x)=P_1\tilde\Gamma^au(t,x)$.
\begin{equation}\label{sharp:decay4}
\Big\|\chi(\frac{|x|}{\w{t}})\nabla^2P_1\tilde\Gamma^au(t,x)\Big\|_{L^2}
\ls\Big\|\nabla\Big(\chi(\frac{|x|}{\w{t}})
P_1\nabla\tilde\Gamma^au(t,x)\Big)\Big\|_{L^2}
+\w{t}^{-1}\|\nabla P_1\tilde\Gamma^au\|_{L^2(|x|\le\w{t}/2)}.
\end{equation}
Applying the first inequality in \eqref{div:curl:ineq} to $\ds\chi(\frac{|x|}{\w{t}})P_1\nabla\tilde\Gamma^au$ yields
\begin{equation}\label{sharp:decay5}
\Big\|\nabla\Big(\chi(\frac{|x|}{\w{t}})P_1\nabla\tilde\Gamma^au(t,x)\Big)\Big\|_{L^2}
\ls\Big\|\chi(\frac{|x|}{\w{t}})\nabla\dive\tilde\Gamma^au(t,x)\Big\|_{L^2}
+\w{t}^{-1}\|\nabla P_1\tilde\Gamma^au\|_{L^2(|x|\le\w{t}/2)}.
\end{equation}
Substituting \eqref{sharp:decay3} and \eqref{sharp:decay5} into \eqref{sharp:decay4} derives
\begin{equation}\label{sharp:decay6}
\w{t}^2\Big\|\chi(\frac{|x|}{\w{t}})\nabla^2P_1\tilde\Gamma^au(t,x)\Big\|_{L^2}
\ls E_{|a|+2}(t)+\cX_{|a|+2}(t)+\cY_{|a|+2}(t)+\w{t}\cW_{|a|}(t).
\end{equation}
Then \eqref{sharp:decay1} is achieved by plugging \eqref{sharp:decay3} and \eqref{sharp:decay6}
into \eqref{sharp:decay2} with $U=P_1\tilde\Gamma^au$ or $U=\Gamma^a\sigma$.
\end{proof}

\begin{lemma}\label{lem:pw:P2u}
For the multi-index $a$ with $|a|\le N-2$, it holds that
\begin{equation}\label{pointwise:P2u}
\|P_2\tilde\Gamma^au(t,x)\|_{L_x^\infty}\ls\cW_{|a|+1}(t).
\end{equation}
\end{lemma}
\begin{proof}
It concludes from the standard Sobolev embedding theorems $W^{1,6}(\R^3)\hookrightarrow L^\infty(\R^3)$,
$\dot H^1(\R^3)\hookrightarrow L^6(\R^3)$ and \eqref{div:curl:ineq} that
\begin{equation*}
\begin{split}
\|P_2U(t,x)\|_{L_x^\infty}
&\ls\|P_2U(t,x)\|_{L_x^6}+\|\nabla P_2U(t,x)\|_{L_x^6}\\
&\ls\|\nabla P_2U(t,x)\|_{L_x^2}+\|\nabla^2P_2U(t,x)\|_{L_x^2}\\
&\ls\|\curl U(t,x)\|_{L_x^2}+\|\nabla\curl U(t,x)\|_{L_x^2}.
\end{split}
\end{equation*}
Thus, \eqref{pointwise:P2u} is proved by choosing $U=\tilde\Gamma^au$.
\end{proof}
Combining Lemma~\ref{lem:sharp:decay} and \ref{lem:pw:P2u} implies the following corollary.
\begin{corollary}\label{coro:sharp:decay}
For $|a|\le N-2$ and $|x|\le\w{t}/4$, it holds that
\begin{equation}\label{sharp:decay7}
|\tilde\Gamma^au(t,x)|+|\Gamma^a\sigma(t,x)|
\ls\w{t}^{-\frac32}[E_{|a|+2}(t)+\cX_{|a|+2}(t)+\cY_{|a|+2}(t)]+\cW_{|a|+1}(t).
\end{equation}
\end{corollary}

\section{Estimates of the auxiliary energies and the good unknown $g$}\label{sect3}

\subsection{Estimates of the auxiliary energies $\cX_N(t)$ and $\cY_N(t)$}\label{sect:aux:energy}
\begin{lemma}[Weighted $\dot H_x^1$ estimate]\label{lem:H1norm}
Under bootstrap assumptions \eqref{bootstrap}, for the integer $m$ with $1\le m\le N$, it holds that
\begin{equation}\label{H1norm}
\cX_m(t)\ls E_m(t)+\cW_{m-1}(t).
\end{equation}
\end{lemma}
\begin{proof}
For the multi-index $|a|\le m-1$, it follows from the equations in \eqref{high:eqn} and direct computations that
\begin{equation}\label{weighted:identity1}
\begin{split}
(|x|^2-t^2)\p_t\tilde\Gamma^au_i&=|x|^2(\cQ^a_{2i}-\p_i\Gamma^a\sigma)
-t\cS\tilde\Gamma^au_i+tx_j\p_j\tilde\Gamma^au_i\\
&=|x|^2\cQ^a_{2i}-x_j(x_j\p_i-x_i\p_j)\Gamma^a\sigma-x_i\cS\Gamma^a\sigma
+tx_i\p_t\Gamma^a\sigma-t\cS\tilde\Gamma^au_i\\
&\quad+tx_j(\p_j\tilde\Gamma^au_i-\p_i\tilde\Gamma^au_j)
+t(x_j\p_i-x_i\p_j)\tilde\Gamma^au_j+tx_i\dive\tilde\Gamma^au\\
&=|x|^2\cQ^a_{2i}-x_j\Omega_{ji}\Gamma^a\sigma-x_i\cS\Gamma^a\sigma
+tx_i\cQ_1^a-t\cS\tilde\Gamma^au_i\\
&\quad+tx_j\eps_{jik}\curl\tilde\Gamma^au_k+t\Omega_{ji}(\tilde\Gamma^au_j),
\end{split}
\end{equation}
where we have used the Einstein summation, and the fact of $\p_jU_i-\p_iU_j=\eps_{jik}\curl U_k$
with  the volume form $\eps_{ijk}$ being the sign of the arrangement $\{ijk\}$.
Note that the main difference between \eqref{weighted:identity1} and similar equality of $\p_tP_1\tilde\Gamma^au_i$
in \cite{Sideris97} lies in the presence of the vorticity $\curl\tilde\Gamma^au$ in \eqref{weighted:identity1}.

Analogously, we can get
\begin{equation}\label{weighted:identity2}
\begin{split}
(|x|^2-t^2)\p_t\Gamma^a\sigma&=|x|^2\cQ^a_1-x_j\Omega_{ji}(\tilde\Gamma^au_i)
-x_i\cS\tilde\Gamma^au_i-t\cS\Gamma^a\sigma+tx_i\cQ^a_{2i},\\
(|x|^2-t^2)\p_i\Gamma^a\sigma&=x_j\Omega_{ji}\Gamma^a\sigma+x_i\cS\Gamma^a\sigma
-tx_i\cQ^a_1-t^2\cQ^a_{2i}+t\cS\tilde\Gamma^au_i\\
&\quad-tx_j\eps_{jik}\curl\tilde\Gamma^au_k-t\Omega_{ji}(\tilde\Gamma^au_j),\\
(|x|^2-t^2)\dive\tilde\Gamma^au&=x_j\Omega_{ji}(\tilde\Gamma^au_i)
+x_i\cS\tilde\Gamma^au_i-tx_i\cQ^a_{2i}-t^2\cQ^a_1+t\cS\Gamma^a\sigma.
\end{split}
\end{equation}
Due to $\w{|x|-t}\ls1+||x|-t|$, by dividing $|x|+t$ and then taking $L_x^2$ norm on the both sides of \eqref{weighted:identity1} and \eqref{weighted:identity2}, we arrive at
\begin{equation}\label{H1norm1}
\cX_m(t)\ls E_m(t)+\cW_{m-1}(t)+\sum_{|b|+|c|\le m-1}
\|\w{|x|+t}(|Q_1^{bc}|+|Q_2^{bc}|)\|_{L_x^2},
\end{equation}
where $Q_1^{bc},Q_2^{bc}$ are defined in \eqref{Qbc:def}.

Next, we estimate the term $\|\w{|x|+t}(|Q_1^{bc}|+|Q_2^{bc}|)\|_{L_x^2}$ on the right hand side of \eqref{H1norm1}.

Since $|b|+|c|\le N-1\le2N_1-3$, then $|b|\le N_1-1$ or $|c|\le N_1-2$ holds.
For the case of $|c|\le N_1-2$, applying \eqref{pointwise:wave} to $\nabla\Gamma^c\sigma,\nabla\tilde\Gamma^cu$ directly yields
\begin{equation}\label{H1norm2}
\begin{split}
&\sum_{\substack{|b|+|c|\le m-1,\\|c|\le N_1-2}}
\|\w{|x|+t}(|Q_1^{bc}|+|Q_2^{bc}|)\|_{L^2}\\
&\ls E_m(t)[E_{N_1+1}(t)+\cX_{N_1+1}(t)+\w{t}\cW_{N_1}(t)]
\ls E_m(t),
\end{split}
\end{equation}
where we have used assumptions \eqref{bootstrap}.

When $|b|\le N_1-1$, the related integral domain will be divided into two parts of $|x|\ge\w{t}/8$ and $|x|\le\w{t}/8$.

In the region of $|x|\ge\w{t}/8$, by using \eqref{pointwise:wave} to $\Gamma^b\sigma,\tilde\Gamma^bu$, we obtain
\begin{equation}\label{H1norm3}
\begin{split}
&\sum_{\substack{|b|+|c|\le m-1,\\|b|\le N_1-1}}
\|\w{|x|+t}(|Q_1^{bc}|+|Q_2^{bc}|)\|_{L^2(|x|\ge\w{t}/8)}\\
&\ls E_m(t)[E_{N_1+1}(t)+\cX_{N_1+1}(t)+\w{t}\cW_{N_1}(t)]\ls E_m(t).
\end{split}
\end{equation}

In the region of $|x|\le\w{t}/8$, it concludes from \eqref{div:curl:ineq} and \eqref{sharp:decay7} that
\begin{equation}\label{H1norm4}
\begin{split}
&\;\quad\|\w{|x|+t}(|Q_1^{bc}|+|Q_2^{bc}|)\|_{L^2(|x|\le\w{t}/8)}\\
&\ls\|\tilde\Gamma^bu,\Gamma^b\sigma\|_{L^\infty(|x|\le\w{t}/8)}
\Big\{\|\w{|x|-t}\nabla\Gamma^c\sigma\|_{L^2}
+\Big\|\w{|x|-t}\nabla\Gamma^cu\Big\|_{L^2}\Big\}\\
&\ls\Big\{\w{t}^{-\frac32}[E_{|b|+2}(t)+\cX_{|b|+2}(t)+\cY_{|b|+2}(t)]
+\cW_{|b|+1}(t)\Big\}[E_m(t)+\cX_m(t)+\w{t}\cW_{m-1}(t)].
\end{split}
\end{equation}
Substituting the assumptions \eqref{bootstrap} into \eqref{H1norm4} yields
\begin{equation}\label{H1norm5}
\sum_{\substack{|b|+|c|\le m-1,\\|b|\le N_1-1}}
\|\w{|x|+t}(|Q_1^{bc}|+|Q_2^{bc}|)\|_{L^2(|x|\le\w{t}/8)}
\ls E_m(t)+M\ve\cX_m(t)+\cW_{m-1}(t).
\end{equation}
Collecting \eqref{H1norm1}--\eqref{H1norm3}, \eqref{H1norm5} with the smallness of $M\ve$, we have achieved \eqref{H1norm}.
\end{proof}

\begin{lemma}[Weighted $\dot H_x^2$ estimate  away from the conic surface $|x|=\w{t}$]\label{lem:H2norm}
Under bootstrap assumptions \eqref{bootstrap}, for each  integer $m$ with $2\le m\le N$, it holds that
\begin{equation}\label{H2norm}
\cY_m(t)\ls E_m(t)+\w{t}\cW_{m-1}(t).
\end{equation}

\end{lemma}
\begin{proof}
For $|a|\le m-2$, replacing $\tilde\Gamma^a,\Gamma^a,\cQ^a_1,\cQ^a_2$ by $\p_x\tilde\Gamma^a,\p_x\Gamma^a,\p_x\cQ^a_1,\p_x\cQ^a_2$
in the last two identities of \eqref{weighted:identity2}, respectively, one has that
\begin{equation}\label{weighted:identity3}
\begin{split}
(|x|^2-t^2)\p_i\p_x\Gamma^a\sigma&=x_j\Omega_{ji}\p_x\Gamma^a\sigma
+x_i\cS\p_x\Gamma^a\sigma-tx_i\p_x\cQ^a_1-t^2\p_x\cQ^a_{2i}\\
&\quad+t\cS\p_x\tilde\Gamma^au_i-tx_j\eps_{jik}\curl\p_x\tilde\Gamma^au_k
-t\Omega_{ji}(\p_x\tilde\Gamma^au_j),\\
(|x|^2-t^2)\dive\p_x\tilde\Gamma^au&=x_j\Omega_{ji}(\p_x\tilde\Gamma^au_i)
+x_i\cS\p_x\tilde\Gamma^au_i-tx_i\p_x\cQ^a_{2i}-t^2\p_x\cQ^a_1\\
&\quad+t\cS\p_x\Gamma^a\sigma.
\end{split}
\end{equation}
Taking $L^2(|x|\le\w{t}/2)$ norm on the both sides of \eqref{weighted:identity3} yields
\begin{equation}\label{H2norm1}
\begin{split}
\cY_m(t)&\ls\cX_m(t)+\w{t}\cW_{m-1}(t)
+\sum_{|a|\le m-2}\|\w{|x|-t}\nabla\Gamma\tilde\Gamma^au\|_{L^2}\\
&\quad+\sum_{|a|\le m-2}\|\w{t}^2(|\p_x\cQ^a_1|+|\p_x\cQ^a_2|)\|_{L^2(|x|\le\w{t}/2)}\\
&\ls E_m(t)+\cX_m(t)+\w{t}\cW_{m-1}(t)+\sum_{|a|\le m-2}
\|\w{t}^2(|\p_x\cQ^a_1|+|\p_x\cQ^a_2|)\|_{L^2(|x|\le\w{t}/2)},
\end{split}
\end{equation}
where we have used the inequality \eqref{div:curl:ineq}.
By the definitions of $\cQ^a_1,\cQ^a_2$ in \eqref{high:eqn}, we have
\begin{equation}\label{H2norm2}
\begin{split}
\p_x\cQ^a_1&=-\sum_{b+c\le a}C^a_{bc}(\lambda\Gamma^b\sigma\p_x\dive\tilde\Gamma^cu
+\tilde\Gamma^bu\cdot\nabla\p_x\Gamma^c\sigma
+\lambda\p_x\Gamma^b\sigma\dive\tilde\Gamma^cu
+\p_x\tilde\Gamma^bu\cdot\nabla\Gamma^c\sigma),\\
\p_x\cQ^a_2&=-\sum_{b+c\le a}C^a_{bc}(\lambda\Gamma^b\sigma\nabla\p_x\Gamma^c\sigma
+\tilde\Gamma^bu\cdot\nabla\p_x\tilde\Gamma^cu
+\lambda\p_x\Gamma^b\sigma\nabla\Gamma^c\sigma
+\p_x\tilde\Gamma^bu\cdot\nabla\tilde\Gamma^cu).
\end{split}
\end{equation}
It suffices to deal with $\p_x\tilde\Gamma^bu\cdot\nabla\tilde\Gamma^cu$ and $\tilde\Gamma^bu\cdot\nabla\p_x\tilde\Gamma^cu$
in \eqref{H2norm2} since the treatments on the other left terms are similar.

In view of $|b|+|c|\le N-2\le2N_1-4$, then $|b|\le N_1-2$ or $|c|\le N_1-2$ holds.
Using \eqref{pointwise:wave} to $\p_x\tilde\Gamma^bu$ or $\nabla\tilde\Gamma^cu$ directly leads to
\begin{equation}\label{H2norm3}
\begin{split}
&\;\quad\w{t}\sum_{b+c=a}\|\p_x\tilde\Gamma^bu\cdot\nabla
\tilde\Gamma^cu\|_{L^2(|x|\le\w{t}/2)}\\
&\ls\sum_{|b|\le N_1-2}\|\p_x\tilde\Gamma^bu\|_{L^\infty}
\|\w{|x|-t}\nabla\tilde\Gamma^cu\|_{L^2}
+\sum_{|c|\le N_1-2}\|\w{|x|-t}\p_x\tilde\Gamma^bu\|_{L^2}
\|\nabla\tilde\Gamma^cu\|_{L^\infty}\\
&\ls\w{t}^{-1}[E_{m-1}(t)+\cX_{m-1}(t)+\w{t}\cW_{m-1}(t)]
[E_{N_1+1}(t)+\cX_{N_1+1}(t)+\w{t}\cW_{N_1}(t)]\\
&\ls\w{t}^{-1}[E_{m-1}(t)+\cX_{m-1}(t)+\w{t}\cW_{m-1}(t)],
\end{split}
\end{equation}
where we have used assumptions \eqref{bootstrap} and inequality \eqref{div:curl:ineq}.

To treat $\tilde\Gamma^bu\cdot\nabla\p_x\tilde\Gamma^cu$, we divide the integral region $|x|\le\w{t}/2$ into
two parts of $\w{t}/4\le|x|\le\w{t}/2$ and $|x|\le\w{t}/4$.

In the region $\w{t}/4\le|x|\le\w{t}/2$, it follows from \eqref{bootstrap}, \eqref{div:curl:ineq} and \eqref{pointwise:wave} that
\begin{equation}\label{H2norm4}
\begin{split}
&\;\quad\|\tilde\Gamma^bu\cdot\nabla
\p_x\tilde\Gamma^cu\|_{L^2(\frac{\w{t}}{4}\le|x|\le\frac{\w{t}}{2})}\\
&\ls\w{t}^{-1}\|\tilde\Gamma^bu\|_{L^\infty(\frac{\w{t}}{4}\le|x|\le\frac{\w{t}}{2})}
\|\w{|x|-t}\nabla\p_x\tilde\Gamma^cu\|_{L^2}\\
&\ls\w{t}^{-\frac52}[E_m(t)+\cX_m(t)+\w{t}\cW_{m-1}(t)]
[E_{|b|+2}(t)+\cX_{|b|+2}(t)+\w{t}\cW_{|b|+1}(t)]\\
&\ls\w{t}^{-2}[E_m(t)+\cX_m(t)+\w{t}\cW_{m-1}(t)].
\end{split}
\end{equation}

In the region $|x|\le\w{t}/4$, it concludes from \eqref{div:curl:ineq}, \eqref{sharp:decay6} and \eqref{sharp:decay7} that
\begin{equation}\label{H2norm5}
\begin{split}
&\;\quad\w{t}^2\|\tilde\Gamma^bu\cdot\nabla
\p_x\tilde\Gamma^cu\|_{L^2(|x|\le\w{t}/4)}\\
&\ls\w{t}^2\|\tilde\Gamma^bu\|_{L^\infty(|x|\le\w{t}/4)}\Big\{\|\nabla\p_xP_2\tilde\Gamma^cu\|_{L^2}
+\|\chi(\frac{|x|}{\w{t}})\nabla\p_xP_1\tilde\Gamma^cu\|_{L^2}\Big\}\\
&\ls\Big\{\w{t}^{-\frac32}[E_m(t)+\cX_m(t)+\cY_m(t)]+\cW_{|b|+1}(t)\Big\}\\
&\qquad\times\Big\{\w{t}^2\cW_{|c|+1}(t)+E_m(t)+\cX_m(t)+\cY_m(t)\Big\}\\
&\ls\w{t}^2\cW_{|c|+1}(t)\{M\ve\w{t}^{-\frac32}+M\delta\w{t}^{M'\ve-\frac12}
+\cW_{|b|+1}(t)\}+\cX_m(t)+\cW_{m-1}(t)+M\ve\cY_m(t),
\end{split}
\end{equation}
where we have also used assumptions \eqref{bootstrap}.

Since $|b|\le N_1-1$ or $|c|\le N_1-1$ always holds, by using assumptions \eqref{bootstrap} again, we have
\begin{equation}\label{H2norm6}
\w{t}^2\cW_{|b|+1}(t)\cW_{|c|+1}(t)\ls\w{t}\cW_{m-1}.
\end{equation}
Therefore, combining \eqref{H2norm1}--\eqref{H2norm6} with the smallness of $M\ve$ derives \eqref{H2norm}.
\end{proof}

\subsection{Estimates of the good unknown $g$}\label{sect:good}

\begin{lemma}\label{lem:good:L2norm}
Under bootstrap assumptions \eqref{bootstrap}, for the good unknown $g$ defined by \eqref{goodunknown:def} and $m\le N-1$, it holds that
\begin{equation}\label{good:L2norm}
\sum_{|a|\le m}\|\w{|x|}\nabla\tilde\Gamma^ag\|^2_{L^2(|x|\ge\w{t}/8)}
\ls E_{m+1}(t)+\cW_m(t).
\end{equation}
\end{lemma}
\begin{proof}
According to \eqref{radial:angular}, it only needs to deal with $r\p_r\tilde\Gamma^ag_i=x_j\p_j\tilde\Gamma^a(u-\sigma\omega)_i$
in order to estimate $\w{|x|}\nabla\tilde\Gamma^ag$ in the left hand side of \eqref{good:L2norm}.
It is deduced from direct computation that there exist some bounded smooth functions $f_i^{a,b}(x)$
and $f_{ij}^{a,b}(x)$ in $|x|\ge1/8$ such that
\begin{equation}\label{Gamma:sigma:omega}
\begin{split}
\tilde\Gamma^a(\sigma\omega)_i&=\omega_i\Gamma^a\sigma
+\frac{1}{|x|}\sum_{b+c\le a}f_i^{a,b}(x)\Gamma^c\sigma,\\
\p_j\tilde\Gamma^a(\sigma\omega)_i&=\omega_i\p_j\Gamma^a\sigma
+\frac{1}{|x|}\sum_{b+c\le a}[f_i^{a,b}(x)\p_j\Gamma^c\sigma+f_{ij}^{a,b}(x)\Gamma^c\sigma].
\end{split}
\end{equation}
Thereafter, we arrive at
\begin{equation}\label{good:L2norm1}
\begin{split}
&\qquad r\p_r\tilde\Gamma^ag_i+\omega_j\sum_{b+c\le a}
[f_i^{a,b}(x)\p_j\Gamma^c\sigma+f_{ij}^{a,b}(x)\Gamma^c\sigma]\\
&=x_j\p_j\tilde\Gamma^au_i-\omega_ix_j\p_j\Gamma^a\sigma\\
&=x_j(\p_j\tilde\Gamma^au_i-\p_i\tilde\Gamma^au_j)+(x_j\p_i-x_i\p_j)\tilde\Gamma^au_j
+x_i\dive\tilde\Gamma^au-\omega_i\cS\Gamma^a\sigma+\omega_it\p_t\Gamma^a\sigma\\
&=x_j\eps_{jik}\curl\tilde\Gamma^au_k+\Omega_{ji}(\tilde\Gamma^au_j)
+x_i\cQ^a_1+\omega_i(t-|x|)\p_t\Gamma^a\sigma-\omega_i\cS\Gamma^a\sigma.
\end{split}
\end{equation}

By taking the $L^2(|x|\ge\w{t}/8)$ norm on the both sides of \eqref{good:L2norm1} and
then substituting \eqref{H1norm}, \eqref{H1norm2} and \eqref{H1norm3} into the resulted inequality, we have
\begin{equation*}
\begin{split}
&\quad\sum_{|a|\le m}\|\w{|x|}\p_r\tilde\Gamma^ag\|_{L^2(|x|\ge\w{t}/8)}\\
&\ls E_{m+1}(t)+\cX_m(t)+\cW_m(t)
+\sum_{|b|+|c|\le|a|}\|\w{|x|}Q_1^{bc}\|_{L^2(|x|\ge\w{t}/8)}\\
&\ls E_{m+1}(t)+\cW_m(t).
\end{split}
\end{equation*}
This together with \eqref{radial:angular} yields
\begin{equation*}
\sum_{|a|\le m}\|\w{|x|}\nabla\tilde\Gamma^ag\|_{L^2(|x|\ge\w{t}/8)}
\ls\sum_{|a|\le m}\|\w{|x|}\p_r\tilde\Gamma^ag\|_{L^2(|x|\ge\w{t}/8)}+E_{m+1}(t).
\end{equation*}
Thus \eqref{good:L2norm} is proved.
\end{proof}

\begin{lemma}\label{lem:good:pw}
Under bootstrap assumptions \eqref{bootstrap}, for $|a|\le N-2$ and $|x|\ge\w{t}/8$, it holds that
\begin{equation}\label{good:pointwise}
\w{|x|+t}^\frac32|\tilde\Gamma^ag(t,x)|\ls E_{|a|+2}(t)+\cW_{|a|+1}(t).
\end{equation}
\end{lemma}
\begin{proof}
Recall (3.19) of \cite{Sideris00} that
\begin{equation}\label{good:pointwise1}
\begin{split}
&\quad\Big(|x|^4|R(|x|)|^2\int_{\SS^2}|U(t,|x|\omega)|^2d\omega\Big)^\frac14\\
&\ls\Big(\int_{|y|\ge|x|}[|R(|y|)|^2|\p_rU(t,y)|^2+|R'(|y|)|^2|U(t,y)|^2]dy\Big)^\frac14
\Big(\int_{|y|\ge|x|}|\tilde\Omega^{\le1}U(t,y)|^2dy\Big)^\frac14.
\end{split}
\end{equation}
By choosing $R(|x|)=\w{|x|}$ and $U(t,x)=\tilde\Omega^{\le1}\tilde\Gamma^ag(t,x)$ in \eqref{good:pointwise1}, 
we can deduce from $W^{1,4}(\SS^2)\hookrightarrow L^\infty(\SS^2)$ that for $|x|\ge\w{t}/8$,
\begin{equation*}
\begin{split}
\w{|x|+t}^\frac32|\tilde\Gamma^ag(t,x)|&\ls\Big(|x|^4\w{|x|}^2\int_{\SS^2}
|\tilde\Omega^{\le1}\tilde\Gamma^ag(t,|x|\omega)|^4d\omega\Big)^\frac14\\
&\ls\sum_{|b|\le|a|+2}\|\tilde\Gamma^bg(t,y)\|_{L^2(|y|\ge\w{t}/8)}
+\sum_{|b|\le|a|+1}\|\w{|y|}\p_r\tilde\Gamma^bg(t,y)\|_{L^2(|y|\ge\w{t}/8)}\\
&\ls E_{|a|+2}(t)+\cW_{|a|+1}(t),
\end{split}
\end{equation*}
where we have used \eqref{good:L2norm} in the last inequality.
This completes the proof of Lemma~\ref{lem:good:pw}.
\end{proof}

\section{Energy estimates}\label{sect4}
Substituting Lemma \ref{lem:H1norm} and \ref{lem:H2norm} into Lemma \ref{lem:pointwise}, \ref{lem:sharp:decay}
and Corollary \ref{coro:sharp:decay} yields
\begin{lemma}\label{lem:decay:summary}
Under bootstrap assumptions \eqref{bootstrap}, for multi-indices $a,b$ with $|a|\le N-2$ and $|b|\le N-3$, it holds that
\begin{align}
\w{|x|}\w{|x|-t}^\frac12(|\tilde\Gamma^au(t,x)|+|\Gamma^a\sigma(t,x)|)
&\ls E_{|a|+2}(t)+\w{t}\cW_{|a|+1}(t),\label{pointwise:wave1}\\
\w{|x|}\w{|x|-t}(|\nabla\tilde\Gamma^bu(t,x)|+|\nabla\Gamma^b\sigma(t,x)|)
&\ls E_{|b|+3}(t)+\w{t}\cW_{|b|+2}(t),\label{pointwise:wave2}
\end{align}
and for $|x|\le\w{t}/4$,
\begin{align}
|\tilde\Gamma^au(t,x)|+|\Gamma^a\sigma(t,x)|
&\ls\w{t}^{-\frac32}E_{|a|+2}(t)+\cW_{|a|+1}(t),\label{sharp:decay}\\
|P_1\tilde\Gamma^au(t,x)|
&\ls\w{t}^{-\frac32}E_{|a|+2}(t)+\w{t}^{-\frac12}\cW_{|a|+1}(t).\label{sharp:decay'}
\end{align}
\end{lemma}

\subsection{Elementary energy estimates}\label{sect:energy}

\begin{lemma}
Under bootstrap assumptions \eqref{bootstrap}, it holds that
if $1<\gamma<3$,
\begin{equation}\label{energy:polytropic}
E^2_N(t')\ls E^2_N(0)+\int_0^{t'}(M\ve\w{t}^{-1}+M\delta)E^2_N(t)dt;
\end{equation}
if $\gamma=-1$,
\begin{equation}\label{energy:Chaplygin}
E^2_N(t')\ls E^2_N(0)+\int_0^{t'}\Big\{(M\ve\w{t}^{-\frac32}+M\delta)E^2_N(t)
+M^2\delta^2\w{t}^{M'\ve-1}E_N(t)\Big\}dt.
\end{equation}
\end{lemma}
\begin{proof}
For $|a|\le N$, multiplying the first equation by $e^q\Gamma^a\sigma$ and the second equation by
$e^q\tilde\Gamma^au$ in \eqref{high:eqn}, where
the ghost weight $e^q=e^{q(|x|-t)}=e^{\arctan(|x|-t)}$ (such a ghost weight is chosen in \cite{Alinhac01a} to
establish the global small data
solutions of 2D quasilinear wave equations with the first null and the second null conditions),
and subsequently adding them, we have
\begin{equation}\label{energy:identity}
\begin{split}
&\frac12\p_t[e^q(|\Gamma^a\sigma|^2+|\tilde\Gamma^au|^2)]
+\dive[e^q(1+\lambda\sigma)\Gamma^a\sigma\tilde\Gamma^au]
+\frac12\dive[e^qu(|\Gamma^a\sigma|^2+|\tilde\Gamma^au|^2)]\\
&+\frac{e^q}{2\w{|x|-t}^2}\sum_{i=1}^3\Big\{|\tilde\Gamma^au_i-\omega_i\Gamma^a\sigma|^2
-u_i\omega_i(|\Gamma^a\sigma|^2+|\tilde\Gamma^au|^2)
-2\lambda\sigma\omega_i\Gamma^a\sigma\tilde\Gamma^au_i\Big\}\\
&=\frac12e^q(|\Gamma^a\sigma|^2+|\tilde\Gamma^au|^2)\dive u
+\lambda e^q\Gamma^a\sigma\tilde\Gamma^au\cdot\nabla\sigma
+\sum_{\substack{b+c=a,\\c<a}}e^qC^a_{bc}(Q_1^{bc}\Gamma^a\sigma
+Q_2^{bc}\cdot\tilde\Gamma^au).
\end{split}
\end{equation}
Integrating the above equality over $[0,t']\times\R^3$ yields
\begin{equation}\label{energy:ineq1}
\begin{split}
&\quad E^2_{|a|}(t')+\sum_{i=1}^3\int_0^{t'}\int\frac{1}{\w{|x|-t}^2}
|\tilde\Gamma^au_i-\omega_i\Gamma^a\sigma|^2dxdt\\
&\ls E^2_{|a|}(0)+\int_0^{t'}\int\Big\{|I^a|+\sum_{\substack{b+c=a,\\c<a}}
(|Q_1^{bc}\Gamma^a\sigma|+|Q_2^{bc}\cdot\tilde\Gamma^au|)\Big\}dxdt,
\end{split}
\end{equation}
where
\begin{equation}\label{Ia:def}
\begin{split}
I^a:=&~(|\Gamma^a\sigma|^2+|\tilde\Gamma^au|^2)\dive u
+2\lambda\Gamma^a\sigma\tilde\Gamma^au\cdot\nabla\sigma\\
&+\frac{1}{\w{|x|-t}^2}\sum_{i=1}^3
\Big\{u_i\omega_i(|\Gamma^a\sigma|^2+|\tilde\Gamma^au|^2)
+2\lambda\sigma\omega_i\Gamma^a\sigma\tilde\Gamma^au_i\Big\}.
\end{split}
\end{equation}
To treat the integral of \eqref{Ia:def}, the related integral domain will be divided into two parts of $|x|\le\w{t}/8$ and $|x|\ge\w{t}/8$.

Since $|b|+|c|\le N\le2N_1-2$, then $|b|\le N_1-1$ or $|c|\le N_1-2$ holds.
In the region $|x|\le\w{t}/8$, by using \eqref{sharp:decay} to $I^a,Q_1^{bc},Q_2^{bc}$ with assumptions \eqref{bootstrap}, we obtain
\begin{equation}\label{Ia:ineq1}
\int_{|x|\le\w{t}/8}|I^a|dx\ls E^2_N(t)\{M\ve\w{t}^{-\frac32}+M\delta\}
\end{equation}
and
\begin{equation}\label{Qbc:ineq1}
\begin{split}
&\quad\int_{|x|\le\w{t}/8}(|Q_1^{bc}\Gamma^a\sigma|+|Q_2^{bc}\cdot\tilde\Gamma^au|)dx\\
&\ls E^2_N(t)\Big\{\sum_{|b|\le N_1-1}[\w{t}^{-\frac32}E_{|b|+2}(t)+\cW_{|b|+1}(t)]\\
&\qquad+\sum_{|c|\le N_1-2}[\w{t}^{-\frac32}E_{|c|+3}(t)+\cW_{|c|+2}(t)]\Big\}\\
&\ls E^2_N(t)\{M\ve\w{t}^{-\frac32}+M\delta\}.
\end{split}
\end{equation}

For the polytropic gases of $1<\gamma<3$, applying \eqref{pointwise:wave1} and \eqref{pointwise:wave2} to $I^a,Q_1^{bc},Q_2^{bc}$
in the region $|x|\ge\w{t}/8$ yields
\begin{equation}\label{Qbc:ineq2}
\int_{|x|\ge\w{t}/8}\Big\{|I^a|+\sum_{\substack{b+c=a,\\c<a}}
(|Q_1^{bc}\Gamma^a\sigma|+|Q_2^{bc}\cdot\tilde\Gamma^au|)\Big\}dx
\ls E^2_N(t)\{M\ve\w{t}^{-1}+M\delta\}.
\end{equation}
Substituting \eqref{Ia:ineq1}--\eqref{Qbc:ineq2} into \eqref{energy:ineq1} implies \eqref{energy:polytropic}.

Next, we turn to the proof of \eqref{energy:Chaplygin}.
In this case, $\gamma=-1$ and $\lambda=\frac{\gamma-1}{2}=-1$ hold.
We point out that in the region $|x|\ge\w{t}/8$, the null condition structures of nonlinearities and
the estimates in subsection 3.2 will play a crucial role.
According to the definition of good unknown \eqref{goodunknown:def} and identities \eqref{Gamma:sigma:omega},
we easily get that
\begin{equation}\label{null:structure}
\begin{split}
\tilde\Gamma^bu_i&=\tilde\Gamma^bg_i+\omega_i\Gamma^b\sigma
+\frac{1}{|x|}\sum_{b_1+b_2\le b}f_i^{b,b_1}(x)\Gamma^{b_2}\sigma,\\
\p_j\tilde\Gamma^cu_i&=\p_j\tilde\Gamma^cg_i+\omega_i\p_j\Gamma^c\sigma
+\frac{1}{|x|}\sum_{c_1+c_2\le c}[f_i^{c,c_1}(x)\p_j\Gamma^{c_2}\sigma
+f_{ij}^{c,c_1}(x)\Gamma^{c_2}\sigma].
\end{split}
\end{equation}
At first, we deal with $I^a$ defined by \eqref{Ia:def}.
It concludes from the definition of good unknown $g$ that
\begin{equation}\label{Ia:ineq2}
\begin{split}
I^a=&\sum_{i=1}^3\Big\{\dive u|\tilde\Gamma^au_i-\omega_i\Gamma^a\sigma|^2
+2\Gamma^a\sigma\tilde\Gamma^au_i(\omega_i\dive u-\p_i\sigma)\Big\}\\
&+\frac{1}{\w{|x|-t}^2}\sum_{i,j=1}^3
\Big\{u_i\omega_i|\tilde\Gamma^au_j-\omega_j\Gamma^a\sigma|^2
+2g_i\omega_i\omega_j\Gamma^a\sigma\tilde\Gamma^au_j\Big\}.
\end{split}
\end{equation}
It follows from \eqref{radial:angular} and the second equality of \eqref{null:structure} that
\begin{equation}\label{Ia:ineq3}
\omega_i\dive u-\p_i\sigma=\omega_i\Big\{\p_jg_j
+\frac{1}{|x|}[f_j^{0,0}(x)\p_j\sigma+f_{jj}^{0,0}(x)\sigma]\Big\}+\frac{1}{|x|}\Omega\sigma.
\end{equation}
Applying \eqref{good:pointwise}, \eqref{pointwise:wave1} and \eqref{pointwise:wave2} to \eqref{Ia:ineq2}
and \eqref{Ia:ineq3} derives
\begin{equation}\label{Ia:ineq4}
\int_{|x|\ge\w{t}/8}|I^a|dx\ls\sum_{i=1}^3\int\frac{M\ve}{\w{|x|-t}^2}
|\tilde\Gamma^au_i-\omega_i\Gamma^a\sigma|^2dx
+E^2_N(t)\{M\ve\w{t}^{-\frac32}+M\delta\},
\end{equation}
where we have also used the Young's inequality.

Next, we focus on the treatments of $Q_1^{bc}$ and $Q_2^{bc}$ defined by \eqref{Qbc:def} with $\gamma=-1$.

\noindent For $|c|\le N_1-2$, by using the second equality of \eqref{null:structure}, we arrive at
\begin{equation}\label{Qbc:ineq3}
\begin{split}
Q_{2i}^{bc}&=-\tilde\Gamma^bu_j\p_j\tilde\Gamma^cu_i
+\Gamma^b\sigma\p_i\Gamma^c\sigma\\
&=-\p_j\tilde\Gamma^cu_i(\tilde\Gamma^bu_j-\omega_j\Gamma^b\sigma)
+\Gamma^b\sigma(\p_i\Gamma^c\sigma-\omega_j\p_j\tilde\Gamma^cu_i)\\
&=-\p_j\tilde\Gamma^cu_i(\tilde\Gamma^bu_j-\omega_j\Gamma^b\sigma)
+\frac{1}{|x|}\Gamma^b\sigma\Omega\Gamma^c\sigma\\
&\quad-\omega_j\Gamma^b\sigma\Big\{\p_j\tilde\Gamma^cg_i
+\frac{1}{|x|}\sum_{c_1+c_2\le c}[f_i^{c,c_1}(x)\p_j\Gamma^{c_2}\sigma
+f_{ij}^{c,c_1}(x)\Gamma^{c_2}\sigma]\Big\},
\end{split}
\end{equation}
and
\begin{equation}\label{Qbc:ineq4}
\begin{split}
Q_1^{bc}&=\Gamma^b\sigma\Big\{\p_i\tilde\Gamma^cg_i
+\frac{1}{|x|}\sum_{c_1+c_2\le c}[f_i^{c,c_1}(x)\p_i\Gamma^{c_2}\sigma
+f_{ii}^{c,c_1}(x)\Gamma^{c_2}\sigma]\Big\}\\
&\quad-\p_i\Gamma^c\sigma(\tilde\Gamma^bu_i-\omega_i\Gamma^b\sigma).
\end{split}
\end{equation}
Applying \eqref{pointwise:wave2} to $\nabla\tilde\Gamma^cu,\nabla\Gamma^c\sigma$ and then using the Young's inequality
to the resulted inequality yield
\begin{equation}\label{Qbc:ineq5}
\begin{split}
&\sum_{|c|\le N_1-2}\int_{|x|\ge\w{t}/8}
|\nabla\tilde\Gamma^cu(\tilde\Gamma^bu_j-\omega_j\Gamma^b\sigma)|dx\\
&\ls\int\frac{M\ve}{\w{|x|-t}^2}|\tilde\Gamma^bu_j-\omega_j\Gamma^b\sigma|^2dx
+E^2_N(t)\{M\ve\w{t}^{-2}+M\delta\}.
\end{split}
\end{equation}
Similarly to the proof of \eqref{Qbc:ineq5}, it concludes from \eqref{good:pointwise}, \eqref{pointwise:wave1}
and \eqref{pointwise:wave2} with \eqref{Qbc:ineq3} and \eqref{Qbc:ineq4} that
\begin{equation}\label{Qbc:ineq6}
\begin{split}
&\sum_{\substack{b+c=a,\\|c|\le N_1-2}}\int_{|x|\ge\w{t}/8}
(|Q_1^{bc}\Gamma^a\sigma|+|Q_2^{bc}\cdot\tilde\Gamma^au|)dx\\
&\ls\sum_{b\le a}\sum_{i=1}^3\int\frac{M\ve}{\w{|x|-t}^2}
|\tilde\Gamma^bu_i-\omega_i\Gamma^b\sigma|^2dx
+E^2_N(t)\{M\ve\w{t}^{-\frac32}+M\delta\}.
\end{split}
\end{equation}
For $|b|\le N_1-1$, substituting the first equality of \eqref{null:structure} into \eqref{Qbc:ineq3}
and \eqref{Qbc:ineq4} yields
\begin{equation}\label{Qbc:ineq7}
\begin{split}
Q_{2i}^{bc}&=-\p_j\tilde\Gamma^cu_i\Big\{\tilde\Gamma^bg_j
+\frac{1}{|x|}\sum_{b_1+b_2\le b}f_j^{b,b_1}(x)\Gamma^{b_2}\sigma\Big\}
+\frac{1}{|x|}\Gamma^b\sigma\Omega\Gamma^c\sigma\\
&\quad-\omega_j\Gamma^b\sigma\Big\{\p_j\tilde\Gamma^cg_i
+\frac{1}{|x|}\sum_{c_1+c_2\le c}[f_i^{c,c_1}(x)\p_j\Gamma^{c_2}\sigma
+f_{ij}^{c,c_1}(x)\Gamma^{c_2}\sigma]\Big\},
\end{split}
\end{equation}
and
\begin{equation}\label{Qbc:ineq8}
\begin{split}
Q_1^{bc}&=\Gamma^b\sigma\Big\{\p_i\tilde\Gamma^cg_i
+\frac{1}{|x|}\sum_{c_1+c_2\le c}[f_i^{c,c_1}(x)\p_i\Gamma^{c_2}\sigma
+f_{ii}^{c,c_1}(x)\Gamma^{c_2}\sigma]\Big\}\\
&\quad-\p_i\Gamma^c\sigma\Big\{\tilde\Gamma^bg_i
+\frac{1}{|x|}\sum_{b_1+b_2\le b}f_i^{b,b_1}(x)\Gamma^{b_2}\sigma\Big\}.
\end{split}
\end{equation}
It is deduced from \eqref{good:L2norm}, \eqref{good:pointwise}, \eqref{pointwise:wave1} and \eqref{pointwise:wave2} that
\begin{equation}\label{Qbc:ineq9}
\begin{split}
&\sum_{\substack{b+c=a,\\c<a,|b|\le N_1-1}}\int_{|x|\ge\w{t}/8}
(|Q_1^{bc}\Gamma^a\sigma|+|Q_2^{bc}\cdot\tilde\Gamma^au|)dx\\
&\ls E^2_N(t)\{M\ve\w{t}^{-\frac32}+M\delta\}
+M^2\delta^2\w{t}^{M'\ve-1}E_N(t).
\end{split}
\end{equation}
Substituting \eqref{Ia:ineq1}, \eqref{Qbc:ineq1}, \eqref{Ia:ineq4}, \eqref{Qbc:ineq6} and \eqref{Qbc:ineq9} into \eqref{energy:ineq1}
derives \eqref{energy:Chaplygin}.
\end{proof}

\subsection{Energy estimates of the vorticity}\label{sect:vorticity}
\begin{lemma}\label{lem:curl:energy}
Under bootstrap assumptions \eqref{bootstrap}, it holds that
\begin{align}
&\cW^2_{N_1}(t')\ls\cW^2_{N_1}(0)
+\int_0^{t'}(M\ve\w{t}^{-\frac43}+M\delta\big)\cW^2_{N_1}(t)dt,\label{curl:low}\\
&\cW^2_{N-1}(t')\ls\cW^2_{N-1}(0)
+\int_0^{t'}(M\ve\w{t}^{-1}+M\delta)\cW^2_{N-1}(t)dt.\label{curl:high}
\end{align}
\end{lemma}

\begin{proof}
It is easy to find the equation of vorticity as follows
\begin{equation}\label{curl:eqn}
(\p_t+u\cdot\nabla)\curl u=\curl u\cdot\nabla u-\curl u\dive u.
\end{equation}
By acting $(\cS+1)^{a_s}\tilde Z^{a_z}$ on the equation \eqref{curl:eqn}, we can find the equation of $\curl\tilde\Gamma^au$:
\begin{equation}\label{curl:eqn:high}
(\p_t+u\cdot\nabla)\curl\tilde\Gamma^au
=\sum_{\substack{b+c=a,\\c<a}}J_1^{bc}+\sum_{b+c=a}J_2^{bc},
\end{equation}
where
\begin{equation}\label{Jbc:def}
\begin{split}
&J_1^{bc}:=\tilde\Gamma^bu\cdot\nabla\curl\tilde\Gamma^cu,\\
&J_2^{bc}:=\curl\tilde\Gamma^cu\cdot\nabla\tilde\Gamma^bu
-\curl\tilde\Gamma^cu\dive\tilde\Gamma^bu.
\end{split}
\end{equation}
Multiplying \eqref{curl:eqn:high} by $\w{|x|}^2e^{q(|x|-t)}\curl\tilde\Gamma^au$ and then integrating the
resulted equality over $[0,t']\times\R^3$ yield
\begin{equation}\label{curl:ineq1}
\begin{split}
&\quad \cW^2_{|a|}(t')+\int_0^{t'}\int
\frac{|\w{|x|}\curl\tilde\Gamma^au|^2}{\w{|x|-t}^2}dxdt\\
&\ls \cW^2_{|a|}(0)+\int_0^{t'}\int|\w{|x|}\curl\tilde\Gamma^au|^2
\Big\{|\dive u|+\frac{|u|}{\w{|x|}}+\frac{|u|}{\w{|x|-t}^2}\Big\}dxdt\\
&\quad+\int_0^{t'}\int\w{|x|}^2|\curl\tilde\Gamma^au|
\Big\{\sum_{\substack{b+c=a,\\c<a}}|J_1^{bc}|+\sum_{b+c=a}|J_2^{bc}|\Big\}dxdt,
\end{split}
\end{equation}
where we have used integration by parts with respect to the space variables.

For the integrand in the second line of \eqref{curl:ineq1}, it concludes from \eqref{pointwise:wave1},
\eqref{pointwise:wave2}, \eqref{sharp:decay} and Young's inequality that
\begin{equation}\label{curl:ineq2}
\begin{split}
|\dive u|+\frac{|u|}{\w{|x|}}+\frac{|u|}{\w{|x|-t}^2}
&\ls M\ve\w{t}^{-\frac32}+\frac{M\ve}{\w{|x|-t}^2}+M\delta,\quad |x|\ge\w{t}/4,\\
|\dive u|+|u|&\ls M\ve\w{t}^{-\frac32}+M\delta,\qquad\qquad\qquad~ |x|\le\w{t}/4.
\end{split}
\end{equation}

Next, we deal with $J_1^{bc}$ and $J_2^{bc}$ in the third line of \eqref{curl:ineq1}.
Similarly to the estimates in the former subsections, the related integral
domain is also divided into two parts of $|x|\ge\w{t}/4$ and $|x|\le\w{t}/4$.

\noindent\underline{$J_1^{bc}$ with $|b|\le|a|-1\le N-2$ in the region $|x|\le\w{t}/4$:}
Due to $|b|+|c|\le N-1\le 2N_1-3$, then $|b|\le N_1-1$ or $|c|\le N_1-1$ holds.
Therefore, it follows from \eqref{sharp:decay} directly that
\begin{equation}\label{J1:ineq1}
\begin{split}
\sum_{c<a,|b|\le|a|-1}\|\w{|x|}J_1^{bc}\|_{L^2(|x|\le\w{t}/4)}
&\ls\sum_{c<a,|b|\le|a|-1}
\cW_{|c|+1}(t)\|\tilde\Gamma^bu\|_{L^\infty(|x|\le\w{t}/4)}\\
&\ls\sum_{c<a,|b|\le|a|-1}\cW_{|c|+1}(t)
\{\w{t}^{-\frac32}E_{|b|+2}(t)+\cW_{|b|+1}(t)\}\\
&\ls\cW_{|a|}(t)\{M\ve\w{t}^{-\frac32}+M\delta\}.
\end{split}
\end{equation}
\noindent\underline{$J_1^{bc}$ with $|b|\le|a|-1$ in the region $|x|\ge\w{t}/4$:}
By using \eqref{pointwise:wave1} to $\tilde\Gamma^bu$ and taking the Young's inequality, 
we easily get that for $|b|\le|a|-1$,
\begin{equation}\label{J1:ineq2}
\begin{split}
|\tilde\Gamma^bu(t,x)|
&\ls\w{t}^{-1}\w{|x|-t}^{-\frac12}\{E_{|b|+2}(t)+\w{t}\cW_{|b|+1}(t)\}\\
&\ls M\ve\w{t}^{-1}\w{|x|-t}^{-\frac12}+\cW_{|b|+1}(t)\\
&\ls M\ve\w{t}^{-\frac43}+M\ve\w{|x|-t}^{-2}+\cW_{|b|+1}(t).
\end{split}
\end{equation}
For $J_2^{bc}$, note that $|b|\le N_1-2$ or $|c|\le N_1$ holds.
Subsequently, similarly to \eqref{J1:ineq1} and \eqref{J1:ineq2}, we have that for $|b|\le|a|-2\le N-3$,
\begin{equation}\label{J2:ineq1}
\begin{split}
\sum_{|b|\le|a|-2}\|\w{|x|}J_2^{bc}\|_{L^2(|x|\le\w{t}/4)}
&\ls\sum_{|b|\le|a|-2}\cW_{|c|}(t)\|\nabla\tilde\Gamma^bu\|_{L^\infty(|x|\le\w{t}/4)}\\
&\ls\sum_{|b|\le|a|-2}\cW_{|c|}(t)\{\w{t}^{-\frac32}E_{|b|+3}(t)+\cW_{|b|+2}(t)\}\\
&\ls\cW_{|a|}(t)\{M\ve\w{t}^{-\frac32}+M\delta\},
\end{split}
\end{equation}
and for $|x|\ge\w{t}/4$,
\begin{equation}\label{J2:ineq2}
\begin{split}
|\nabla\tilde\Gamma^bu(t,x)|
&\ls\w{t}^{-1}\w{|x|-t}^{-1}\{E_{|b|+3}(t)+\w{t}\cW_{|b|+2}(t)\}\\
&\ls M\ve\w{t}^{-2}+M\ve\w{|x|-t}^{-2}+\cW_{|b|+2}(t).
\end{split}
\end{equation}
By \eqref{J1:ineq1}--\eqref{J2:ineq2} and the fact of $|a|\le N_1\le N-3(\Rightarrow N_1\ge5)$,
in order to achieve the lower order energy estimate
\eqref{curl:low}, it remains to control $J_1^{bc}$ with $b=a$ and $J_2^{bc}$ with $|b|\ge|a|-1$ 
in the region $|x|\le\w{t}/4$.
For this purpose, by the Helmholtz decomposition \eqref{Helmholtz}, we have
\begin{equation}\label{J1:low1}
\begin{split}
&\sum_{b=a,c<a}\|\w{|x|}J_1^{bc}\|_{L^2(|x|\le\w{t}/4)}\\
&\ls\sum_{b=a,c<a}\Big\{\cW_{|c|+1}(t)\|P_1\tilde\Gamma^bu\|_{L^\infty(|x|\le\w{t}/4)}\\
&\quad+\|\w{|x|}^2\nabla\curl\tilde\Gamma^cu\|_{L^\infty}
\|\w{|x|}^{-1}P_2\tilde\Gamma^bu\|_{L^2}\Big\}.
\end{split}
\end{equation}
Applying the Hardy's inequality and \eqref{div:curl:ineq} to the last term in \eqref{J1:low1} yield
\begin{equation}\label{J1:low2}
\|\w{|x|}^{-1}P_2\tilde\Gamma^bu\|_{L^2}\ls\|\nabla P_2\tilde\Gamma^bu\|_{L^2}
\ls\|\curl P_2\tilde\Gamma^bu\|_{L^2}\ls\cW_{|a|}(t).
\end{equation}
Thereafter, by plugging \eqref{pointwise:curl}, \eqref{sharp:decay'} and \eqref{J1:low2} into \eqref{J1:low1}, we obtain
\begin{equation}\label{J1:low3}
\begin{split}
&\sum_{b=a,c<a}\|\w{|x|}J_1^{bc}\|_{L^2(|x|\le\w{t}/4)}\\
&\ls\cW_{|a|}(t)\{M\ve\w{t}^{-\frac32}+M\delta\w{t}^{M'\ve-\frac12}\}
+\cW_3(t)\cW_{|a|}(t)\\
&\ls\cW_{N_1}(t)\{M\ve\w{t}^{-\frac32}+M\delta\},
\end{split}
\end{equation}
where we have used the bootstrap assumptions \eqref{bootstrap} with $N_1\ge5$ in the last line of \eqref{J1:low3}.
Analogously, we can get the following estimate of $J_2^{bc}$,
\begin{equation}\label{J2:low1}
\begin{split}
&\sum_{\substack{b+c=a,\\|b|\ge|a|-1}}\|\w{|x|}J_2^{bc}\|_{L^2(|x|\le\w{t}/4)}\\
&\ls\sum_{|b|\le|a|,|c|\le1}\Big\{\cW_{|c|}(t)
\|P_1\nabla\tilde\Gamma^bu\|_{L^\infty(|x|\le\w{t}/4)}
+\|\w{|x|}\curl\tilde\Gamma^cu\|_{L^\infty}\|P_2\nabla\tilde\Gamma^bu\|_{L^2}\Big\}\\
&\ls\cW_{N_1}(t)\{M\ve\w{t}^{-\frac32}+M\delta\}.
\end{split}
\end{equation}
Collecting \eqref{curl:ineq1}--\eqref{J2:ineq2}, \eqref{J1:low3}, \eqref{J2:low1} together with all $|a|\le N_1$,
we eventually achieve
\begin{equation*}
\begin{split}
&\cW^2_{N_1}(t')+\sum_{|a|\le N_1}\int_0^{t'}
\Big\|\frac{\w{|x|}\curl\tilde\Gamma^au}{\w{|x|-t}}\Big\|^2_{L_x^2}dt\\
\ls&~\cW^2_{N_1}(0)+\int_0^{t'}\cW^2_{N_1}(t)\{M\ve\w{t}^{-\frac43}+M\delta\}dt\\
&+M\ve\sum_{|b|\le N_1}\int_0^{t'}
\Big\|\frac{\w{|x|}\curl\tilde\Gamma^bu}{\w{|x|-t}}\Big\|^2_{L_x^2}dt.
\end{split}
\end{equation*}
This together with the smallness of $M\ve$ implies \eqref{curl:low}.

At last, the proof of \eqref{curl:high} is reduced to the case of $N_1+1\le|a|\le N-1$.
In view of \eqref{J1:ineq1}--\eqref{J2:ineq2}, we  only need to treat $J_1^{bc}$ and $J_2^{bc}$ for $|b|\ge|a|-1$.
Note that $|b|\le|a|\le N-1$ and $|c|\le1\le N_1-3\le N-4$ hold.
In the region $|x|\ge\w{t}/4$, applying \eqref{bootstrap} and \eqref{pointwise:curl} directly leads to
\begin{equation}\label{J1:high1}
\begin{split}
\|\w{|x|}J_1^{bc}\|_{L^2(|x|\ge\w{t}/4)}
&\ls\w{t}^{-1}E_{|b|}(t)
\|\w{|x|}^2\nabla\curl\tilde\Gamma^cu\|_{L^\infty}\\
&\ls\w{t}^{-1}E_{|b|}(t)\cW_{|c|+3}(t)\ls M\ve\w{t}^{-1}\cW_{N-1}(t).
\end{split}
\end{equation}
In the region $|x|\le\w{t}/4$, it concludes form the Hardy inequality, \eqref{bootstrap}, \eqref{div:curl:ineq} 
and \eqref{H1norm} that
\begin{equation}\label{J1:high2}
\begin{split}
\|\w{|x|}J_1^{bc}\|_{L^2(|x|\le\w{t}/4)}
&\ls\w{t}^{-1}\Big\|\frac{\w{|x|-t}\tilde\Gamma^bu}{\w{|x|}}\Big\|_{L^2}
\|\w{|x|}^2\nabla\curl\tilde\Gamma^cu\|_{L^\infty}\\
&\ls\w{t}^{-1}\cW_{|c|+3}(t)\{E_{|b|}(t)+\cX_{|b|+1}(t)+\w{t}\cW_{|b|}(t)\}\\
&\ls\w{t}^{-1}\cW_{N-1}(t)E_N(t)+\cW_{N_1}(t)\cW_{N-1}(t)\\
&\ls\cW_{N-1}(t)\{M\ve\w{t}^{-1}+M\delta\}.
\end{split}
\end{equation}
Similarly, we can get the following estimate of $J_2^{bc}$,
\begin{equation}\label{J2:high1}
\begin{split}
\|\w{|x|}J_2^{bc}\|_{L^2}
&\ls\|\w{|x|-t}^{-1}\w{|x|}\curl\tilde\Gamma^cu\|_{L^\infty}
\|\w{|x|-t}\nabla\tilde\Gamma^bu\|_{L^2}\\
&\ls\w{t}^{-1}\cW_{|c|+2}(t)\{E_{|b|+1}(t)+\w{t}\cW_{|b|}(t)\}\\
&\ls\cW_{N-1}(t)\{M\ve\w{t}^{-1}+M\delta\w{t}\}.
\end{split}
\end{equation}
For all $|a|\le N-1$, based on the estimate \eqref{curl:low}, substituting \eqref{curl:ineq2}--\eqref{J2:ineq2} and \eqref{J1:high1}--\eqref{J2:high1} into \eqref{curl:ineq1} yields
\begin{equation*}
\begin{split}
\cW^2_{N-1}(t')+\sum_{|a|\le N-1}\int_0^{t'}
\Big\|\frac{\w{|x|}\curl\tilde\Gamma^au}{\w{|x|-t}}\Big\|^2_{L_x^2}dt\\
\ls\cW^2_{N-1}(0)+\int_0^{t'}\cW^2_{N-1}(t)\big\{M\ve\w{t}^{-1}+M\delta\big\}dt\\
+M\ve\sum_{|b|\le N-1}\int_0^{t'}
\Big\|\frac{\w{|x|}\curl\tilde\Gamma^bu}{\w{|x|-t}}\Big\|^2_{L_x^2}dt.
\end{split}
\end{equation*}
Then \eqref{curl:high} is proved.
\end{proof}

\section{Proof of Theorem ~\ref{thm:main}}\label{sect5}

\begin{proof}[Proof of Theorem ~\ref{thm:main}]

(i) Applying the Growall's inequality to \eqref{energy:polytropic}, \eqref{curl:low}, \eqref{curl:high} and
then combining the resulted inequalities with \eqref{initial}, \eqref{H1norm} and \eqref{H2norm}, we know
that there exist two positive constants $C_1,C_2\ge1$ such that
\begin{equation}\label{polytropic:ineq}
\begin{split}
E_N(t)+\cX_N(t)\le C_1\ve(1+t)^{C_2M\ve},
\qquad&\cY_N(t)\le C_1\ve(1+t)^{C_2M\ve}+C_1\delta(1+t)^{1+C_2M\ve},\\
\cW_{N-1}(t)\le C_1\delta(1+t)^{C_2M\ve},
\qquad&\cW_{N_1}(t)\le C_1\delta.
\end{split}
\end{equation}
Choosing $M=2eC_1$, $M'=2eC_1C_2$, $\kappa_0=\frac{1}{4eC_1C_2}$ and $\ve_0=\delta_0=\frac{1}{4eC_1}$, we then  obtain
that for $t\le T=\min\{e^\frac{\kappa_0}{\ve}-1,\frac{\kappa_0}{\delta}\}$,
\begin{equation*}
\begin{split}
E_N(t)+\cX_N(t)\le\frac12M\ve,
\qquad&\cY_N(t)\le\frac12M\ve+\frac12M\delta(1+t)^{1+M'\ve},\\
\cW_{N-1}(t)\le\frac12M\delta(1+t)^{M'\ve},
\qquad&\cW_{N_1}(t)\le\frac12M\delta.
\end{split}
\end{equation*}
This, together with the local existence of classical solution to \eqref{reducedEuler} (see Chapter 2 of \cite{Majda}),
implies that \eqref{reducedEuler}
with \eqref{polytropic:gas} admits a unique solution $(\sigma,u)\in C([0,T],H^N(\R^3))$.
Hence, the proof of Theorem~\ref{thm:main} (i) is completed.

(ii) Let $\tilde E_N(t):=\sup_{0\le s\le t}E_N(s)$ and for convenience we still denote $\tilde E_N(t)$ as $E_N(t)$.
Then it follows from \eqref{energy:Chaplygin} that
\begin{equation}\label{energy:Chaplygin'}
E_N(t')\ls E_N(0)+\frac{M^2\delta^2\w{t}^{M'\ve}}{M'\ve}
+\int_0^{t'}(M\ve\w{t}^{-\frac32}+M\delta)E_N(t)dt.
\end{equation}
If $\delta\le O(\ve^\frac87)$, then for $t\delta\le\kappa_0$, we have
$\delta\w{t}^{M'\ve}\ls\delta\w{t}^\frac18\ls\delta^\frac78\ls\ve$.
Plugging this inequality into \eqref{energy:Chaplygin'} and utilizing the Growall's inequality to the resulted
inequality with \eqref{curl:low} and \eqref{curl:high}, similarly to the proof of \eqref{polytropic:ineq}, we
know that there exist two positive constants $C_3,C_4\ge1$ such that
\begin{equation*}
\begin{split}
E_N(t)+\cX_N(t)\le C_3\ve(1+\frac{M^2}{M'}),
\qquad&\cY_N(t)\le C_3\ve(1+\frac{M^2}{M'})+C_3\delta(1+t)^{1+C_4M\ve},\\
\cW_{N-1}(t)\le C_3\delta(1+t)^{C_4M\ve},
\qquad&\cW_{N_1}(t)\le C_3\delta.
\end{split}
\end{equation*}
Let $M=4C_3$, $M'=\max\{16C_3^2,4C_3C_4\}$, $\kappa_0=\frac{1}{8C_3}$ and $\ve_0=\delta_0=\frac{1}{8M'}$, 
then for $t\le T=\frac{\kappa_0}{\delta}$, we achieve
\begin{equation}\label{bootstrap'}
\begin{split}
E_N(t)+\cX_N(t)\le\frac12M\ve,
\qquad&\cY_N(t)\le\frac12M\ve+\frac12M\delta(1+t)^{1+M'\ve},\\
\cW_{N-1}(t)\le\frac12M\delta(1+t)^{M'\ve},
\qquad&\cW_{N_1}(t)\le\frac12M\delta.
\end{split}
\end{equation}
If $\delta=O(\ve^{1+\alpha})$ with $0<\alpha<\frac17$ and $M'\ve_0\le\frac{\alpha}{1+\alpha}\le\frac18$,
we find $\delta\w{t}^{M'\ve}\ls\delta^{1-M'\ve_0}\ls\ve$.
Analogously, choosing $M, M', \kappa_0$ as before and $\ve_0=\delta_0=\frac{\alpha}{M'(1+\alpha)}$, 
we can get \eqref{bootstrap'}.
Therefore, \eqref{reducedEuler} with \eqref{Chaplygin:gas} will admit a unique 
solution $(\sigma,u)\in C([0,T],H^N(\R^3))$,
which completes the proof of Theorem~\ref{thm:main} (ii)  by the continuity argument.
\end{proof}


\begin{thebibliography}{99}


\bibitem{Alinhac93} S. Alinhac, {\it Temps de vie des solutions r\'eguli\'eres des \'equations d'Euler
compressibles axisym\'etriques en dimension deux,} Invent. Math. \textbf{111} (1993), 627--670.

\bibitem{Alinhac95} S. Alinhac, {\it   Blowup for nonlinear hyperbolic equations}. Progress in Nonlinear
Differential Equations and their Applications, 17. Birkh\"auser Boston, Inc., Boston, MA, 1995.

\bibitem{Alinhac95-1}  S. Alinhac, {\it  Explosion g\'eom\'etrique pour des syst\'emes
quasi-lin\'eaires}, Amer. J. Math. 117 (1995), no. 4, 987--1017.


\bibitem{Alinhac01a} S. Alinhac, {\it The null condition for quasilinear wave equations in two space
dimensions I,} Invent. Math. \textbf{145} (2001), no. 3, 597--618.

\bibitem{Alinhac01b} S. Alinhac, {\it The null condition for quasilinear wave equations in two space
dimensions II,} Amer. J. Math. \textbf{123} (2001), 1071--1101.

\bibitem{Alinhac10} S. Alinhac, {\it Geometric analysis of hyperbolic differential equations: an introduction.}
 London Mathematical Society Lecture Note Series, 374. Cambridge University Press, Cambridge, 2010. x+118 pp.


\bibitem{BSV-1} T. Buckmaster, S. Shkoller, V. Vicol, {\it Formation of point shocks for 3d compressible Euler}, arXiv:1912.04429 (2019).

\bibitem{BSV-2} T. Buckmaster, S. Shkoller, V. Vicol, {\it Formation of shocks for 2D isentropic compressible Euler}, arXiv:1907.03784 (2019).

\bibitem{BSV-3} T. Buckmaster, S. Shkoller, V. Vicol, {\it Shock formation and vorticity creation for 3d Euler}, arXiv:2006.14789 (2020).

\bibitem{Christodoulou86} D. Christodoulou, {\it Global solutions of nonlinear hyperbolic equations for small initial data}, Comm. Pure Appl. Math. \textbf{39} (1986), no. 2, 267--282.

\bibitem{Christodoulou07} D. Christodoulou, {\it The formation of shocks in 3-dimensional fluids,} EMS Monogr. Math., Eur. Math. Soc., Z\"urich, 2007.

\bibitem{CM14} D. Christodoulou, Miao Shuang, {\it Compressible flow and Euler's equations}, Surveys of Modern Mathematics, \textbf{9}, International Press, Somerville, MA; Higher Education Press, Beijing, 2014.

\bibitem{CF:book} R. Courant, K. O. Friedrichs, {\it Supersonic flow and shock waves,} Interscience Publishers Inc., New York, 1948.


\bibitem{Godin05} P. Godin, {\it The lifespan of a class of smooth spherically symmetric solutions of the compressible Euler equations with variable entropy in three space dimensions,} Arch. Ration. Mech. Anal. \textbf{177} (2005), no. 3, 479--511.

\bibitem{Godin07} P. Godin, {\it Global existence of a class of smooth 3D spherically symmetric flows of Chaplygin gases with variable entropy,} J. Math. Pures Appl. \textbf{87} (2007), 91--117.

\bibitem{GIP16} Guo Yan, A.D. Ionescu, B. Pausader, {\it Global solutions of the Euler-Maxwell two-fluid system in 3D,} Ann. of Math. (2) \textbf{183} (2016), 377--498.

\bibitem{HKSW} G. Holzegel, S. Klainerman, J. Speck, W.W.-Y. Wong, {\it Small-data shock formation in solutions to 3d quasilinear wave equations: An overview,} Journal of Hyperbolic Differential Equations \textbf{13} (2016), no. 01, 1--105.

\bibitem{Hormander97book} L. H\"ormander, {\it Lectures on nonlinear hyperbolic differential equations.} Math\'ematiques \& Applications (Berlin) [Mathematics \& Applications], \textbf{26}. Springer-Verlag, Berlin, 1997. viii+289 pp.

\bibitem{HouYin19} Hou Fei, Yin Huicheng, {\it Global smooth axisymmetric solutions to 2D compressible Euler equations of Chaplygin gases with non-zero vorticity,} J. Differential Equations \textbf{267} (2019), no. 5, 3114--3161.

\bibitem{HouYin20} Hou Fei, Yin Huicheng, {\it On global axisymmetric solutions to 2D compressible full Euler equations of Chaplygin gases,} Discrete Contin. Dyn. Syst. \textbf{40} (2020) no. 3, 1435--1492.

\bibitem{HouYin20jde} Hou Fei, Yin Huicheng, {\it Global small data smooth solutions of 2-D null-form wave equations with non-compactly supported initial data,} J. Differential Equations \textbf{268} (2020), no. 2, 490--512.

\bibitem{HouYin21} Hou Fei, Yin Huicheng, {\it Long time existence of smooth solutions to 2D compressible Euler equations of Chaplygin gases with non-zero vorticity,}  arXiv:2102.12038, Preprint (2021).

\bibitem{IL18} A.D. Ionescu, V. Lie, {\it Long term regularity of the one-fluid Euler-Maxwell system in 3D with vorticity,} Adv. Math. \textbf{325} (2018), 719--769.

\bibitem{John} F. John, {\it Nonlinear wave equations, formation of singularities.}
  Seventh Annual Pitcher Lectures delivered at Lehigh University, Bethlehem, Pennsylvania, April 1989. University Lecture Series, 2. American Mathematical Society, Providence, RI, 1990.



\bibitem{Klainerman} S. Klainerman, {\it The null condition and global existence to nonlinear wave equations,} in: Nonlinear Systems of Partial Differential Equations in Applied Mathematics, Part \textbf{1}, Santa Fe, NM, 1984, in: Lect. Appl. Math., vol. \textbf{23}, Amer. Math. Soc., Providence, RI, 1986, pp. 293--326.


\bibitem{Li} Li Ta-tsien, {\it Global classical solutions for quasilinear hyperbolic systems}, Res. Appl. Math., Vol. \textbf{32}, Wiley/Masson, New York/Paris, 1994.

\bibitem{Lax} P. D.Lax, {\it Hyperbolic systems of conservation laws and the mathematical theory of shock waves,} Conference Board of the Mathematical Sciences Regional Conference Series in Applied Mathematics, No. 11. Society for Industrial and Applied Mathematics, Philadelphia, Pa., 1973.

\bibitem{LukSpeck18} J. Luk, J. Speck, {\it Shock formation in solutions to the 2D compressible Euler equations in the presence of non-zero vorticity,} Invent. Math. \textbf{214} (2018), no. 1, 1--169.

\bibitem{Majda} A. Majda, {\it Compressible fluid flow and systems of conservation laws in several space variables,} Applied Mathematical Sciences, \textbf{53}, Springer-Verlag, New York, 1984.



\bibitem{Secchi04} A. Morando, P. Secchi, {\it On 3D slightly compressible Euler equations},
 Port. Math. (N.S.) 61 (2004), no. 3, 301--316.



\bibitem{Sideris97} T. Sideris, {\it Delayed singularity formation in 2D compressible flow,} Amer. J. Math. \textbf{119} (1997), 371--422.

\bibitem{Sideris00} T. Sideris, {\it Nonresonance and global existence of prestressed nonlinear elastic waves,} Ann. of Math. (2) \textbf{151} (2000), no. 2, 849--874.


\bibitem{Yin} Yin Huicheng, {\it Formation and construction of a shock wave for 3-D compressible Euler equations with the spherical initial data,} Nagoya Math. J., \textbf{175} (2004), 125--164.

\end{thebibliography}
\end{document}